\documentclass[reqno, 10pt]{article}

\usepackage{fullpage}
\usepackage{newlfont}
\usepackage{color}
\usepackage{authblk}

\usepackage{microtype}
\usepackage{graphicx}
\usepackage{subfigure}
\usepackage{booktabs} 
\usepackage{algorithm}
\usepackage{algorithmic}

\usepackage{xr-hyper}

\usepackage{amsmath}
\usepackage{amsfonts}
\usepackage{amssymb}
\usepackage{amsthm}
\usepackage{bbm}

\usepackage{hyperref}

\newtheorem{thm}{Theorem}
\newtheorem{cor}[thm]{Corollary}
\newtheorem{lem}[thm]{Lemma}

\theoremstyle{definition}
\newtheorem{defn}[thm]{Definition}
\newtheorem{assump}[thm]{Assumption}

\theoremstyle{remark}
\newtheorem{rem}[thm]{Remark}

\newcommand{\R} {\mathbb{R}}
\newcommand{\C} {\mathbb{C}}
\newcommand{\N} {\mathbb{N}}

\newcommand{\E} {\mathbb{E}}
\newcommand{\T} {\mathbb{T}}
\newcommand{\p} {\mathbb{P}}

\newcommand{\abs}[1]{\lvert #1 \rvert}
\newcommand{\pbb}[1]{\biggl({#1}\biggr)}
\newcommand{\pb}[1]{\bigl({#1}\bigr)}

\DeclareMathOperator{\Tr}{Tr}

\DeclareMathOperator{\err}{\mathrm{err}}
\DeclareMathOperator{\erfc}{\mathrm{erfc}}

\DeclareMathOperator{\im}{\mathrm{Im}}

\newcommand{\caE}{{\mathcal E}}
\newcommand{\caG}{{\mathcal G}}
\newcommand{\caK}{{\mathcal K}}
\newcommand{\caP}{{\mathcal P}}
\newcommand{\caR}{{\mathcal R}}

\newcommand{\caL}{{\mathcal L}}
\newcommand{\caN}{{\mathcal N}}
\newcommand{\caO}{{\mathcal O}}

\newcommand{\bsa}{{\boldsymbol a}}
\newcommand{\bsb}{{\boldsymbol b}}
\newcommand{\bse}{{\boldsymbol e}}
\newcommand{\bsH}{{\boldsymbol H}}
\newcommand{\bsu}{{\boldsymbol u}}
\newcommand{\bsv}{{\boldsymbol v}}
\newcommand{\bsw}{{\boldsymbol w}}
\newcommand{\bsx}{{\boldsymbol x}}
\newcommand{\bsy}{{\boldsymbol y}}
\newcommand{\bsz}{{\boldsymbol z}}
\newcommand{\bsth}{{\boldsymbol \theta}}

\newcommand{\wt}{\widetilde}

\newcommand{\wh}{\widehat}
\newcommand{\beq}{ \begin{equation} }
\newcommand{\eeq}{ \end{equation} }

\newcommand{\dd}{\mathrm{d}}
\newcommand{\ii}{\mathrm{i}}

\renewcommand{\P}{\mathbb{P}}

\newcommand{\f}{q}
\newcommand{\F}{Q}

\newcommand{\fg}{F_g}
\newcommand{\ef}{E_\f}
\newcommand{\mf}{M_\f}
\newcommand{\vf}{V_\f}
\newcommand{\SNR}{\omega}
\newcommand{\rat}{d_0}

\numberwithin{equation}{section} 
\numberwithin{thm}{section}

\title{Detection of Signal in the Spiked Rectangular Models}

\author{Ji Hyung Jung\footnote{Department of Mathematical Sciences, KAIST, Daejeon, 34141, Korea
		\newline email: \texttt{jhjung66@kaist.ac.kr}},	\;
	Hye Won Chung \footnote{School of Electrical Engineering, KAIST, Daejeon, 34141, Korea
		\newline email: \texttt{hwchung@kaist.ac.kr}}\;
	and Ji Oon Lee\footnote{Department of Mathematical Sciences, KAIST, Daejeon, 34141, and School of Mathematics, KIAS, Seoul, 02455, Korea
		\newline email: \texttt{jioon.lee@kaist.edu}}}

\date{\today}

\begin{document}

\maketitle

\begin{abstract}
We consider the problem of detecting signals in the rank-one signal-plus-noise data matrix models that generalize the spiked Wishart matrices. We show that the principal component analysis can be improved by pre-transforming the matrix entries if the noise is non-Gaussian. As an intermediate step, we prove a sharp phase transition of the largest eigenvalues of spiked rectangular matrices, which extends the Baik--Ben Arous--P\'ech\'e (BBP) transition. We also propose a hypothesis test to detect the presence of signal with low computational complexity, based on the linear spectral statistics, which minimizes the sum of the Type-I and Type-II errors when the noise is Gaussian. 
\end{abstract}

\section{Introduction} \label{sec:intro}

Detecting a low-rank structure or signal in a high-dimensional noisy data is one of the most fundamental problems in statistics and data science \cite{Johnstone2001,Onatski2013,Onatski2014,Abbe2017}.
In the case where the data is a matrix and the signal is a vector, it is natural to consider spiked random matrices, which includes the spiked Wigner matrix and the spiked Wishart matrix. In these models, the signal is in the form of rank-$1$ mean matrix (spiked Wigner matrix) or rank-$1$ perturbation of the identity in the covariance matrix (spiked Wishart matrix). In this paper, we consider the following rectangular random matrix models that generalize the spiked Wishart matrix:

\begin{itemize}
	\item Rectangular matrix \emph{with spiked mean} (additive model): the data matrix is of the form
	\[
	\sqrt{\lambda} \bsu \bsv^T + X,
	\]
	where $X$ is an $M \times N$ random i.i.d. matrix whose entries are centered with variance $N^{-1}$, $\bsu \in \R^M$, $\bsv \in \R^N$ with $\| \bsu \| = \| \bsv \| = 1$. The parameter $\lambda$ corresponds to the signal-to-noise ratio (SNR).
	
	\item Rectangular matrix \emph{with spiked covariance} (multiplicative model): the data matrix is of the form
	\[
	(I+ \lambda \bsu \bsu^T)^{1/2} X,
	\]
	where $X$ is an $M \times N$ random i.i.d. matrix whose entries are centered with variance $N^{-1}$, $\bsu \in \R^M$ with $\| \bsu \| = 1$. The parameter $\lambda$ corresponds to the SNR.
\end{itemize}

Note that in the rectangular matrix with spiked covariance, and also in the rectangular matrix with spiked mean under an additional assumption that the entries of $\bsu$ and $\bsv$ are centered, the population covariance is \[\Sigma=I+ \lambda \bsu \bsu^T.\] In the special case where the entries of $\bsv$ are i.i.d. Gaussians, the two models coincide.

If SNR $\lambda$ is sufficiently large, we can easily detect (and recover) the signal by methods such as principal component analysis (PCA). Even under the high-dimensional assumption $M, N \to \infty$ with $M/N \to \rat \in (0, \infty)$, the signal can be reliably detected by PCA if $\lambda>\sqrt{\rat}$. (For the use of PCA in the high-dimensional setting, we refer to \cite{johnstone2007}.) On the other hand, if $\lambda \in (0, \sqrt{\rat})$, the distribution of the largest eigenvalue coincides with that of the null model $\lambda = 0$. This sharp transition in the behavior of the largest eigenvalue is known as the BBP transition after the seminal work by Baik, Ben Arous, and P\'ech\'e \cite{BBP2005}. (See Section \ref{subsec:BBP}.) 

On the other hand, in the subcritical case $\lambda < \sqrt{\rat}$, if the noise $X$ is Gaussian and the signal $\bsu$ (and also $\bsv$ for a rectangular matrix with spiked mean) is drawn uniformly from the unit sphere, known as the spherical prior, then no test can reliably detect the signal. (See Section \ref{subsec:LR}.) Thus, it is natural to ask the following questions:

\begin{itemize}
	\item Is the threshold for reliable detection (i.e., with probability $1-o(1)$ as $M, N \to \infty$) lower than $\sqrt{\rat}$ if the noise is non-Gaussian?
	
	\item Can we design an efficient algorithm to weakly detect the signal (i.e., better than a random guess) for the subcritical case?
\end{itemize}

We aim to answer these questions in this paper.

\subsection{Main contributions} \label{subsec:main}

Our main contributions are as follows:
\begin{itemize}
	\item We prove that the PCA can be improved by an entrywise transformation if the noise is non-Gaussian, under a mild assumption on the distribution (prior) of the spike.
	
	\item We propose a universal test to detect the presence of signal with low computational complexity, based on the linear spectral statistics (LSS). The test does not require any prior information on the signal, and if the noise is Gaussian the error of the proposed test is optimal.
\end{itemize}

Heuristically, the SNR can be increased through an entrywise transformation and it can be easily seen for a rectangular matrix (additive model) of the form $Y = \sqrt{\lambda} \bsu \bsv^T + X$. If $|u_i v_j| \ll X_{ij}$, then by applying a function $\f$ entrywise to $\sqrt{N} Y$, we obtain a transformed matrix whose entries are
\[ \begin{split}
\f(\sqrt{N} Y_{ij} ) =\f(\sqrt{N} X_{ij} + \sqrt{\lambda N} u_i v_j) \approx \f(\sqrt{N} X_{ij}) + \sqrt{\lambda N} \f'(\sqrt{N} X_{ij}) u_i v_j ,
\end{split} \]
where the approximation is due to the Taylor expansion. It can be shown that the coefficient $\f'(\sqrt{N} X_{ij})$ in the second term in the right side can be replaced by its expectation with negligible error. (See Appendix \ref{sec:proof_trans-mean} for the proof.) Thus,
\[ \begin{split}
\f(\sqrt{N} Y_{ij} ) =\f(\sqrt{N} X_{ij} + \sqrt{\lambda N} u_i v_j) \approx \sqrt{N} \left( \frac{\f(\sqrt{N} X_{ij})}{\sqrt{N}} + \sqrt{\lambda} \E[\f'(\sqrt{N} X_{ij})] u_i v_j \right),
\end{split} \]
and the transformed matrix is of the form $\sqrt{\lambda'} \bsu \bsv^T + \F$ after normalization, which yields another spiked rectangular matrix with different SNR. By optimizing the SNR of the transformed matrix, we find that the SNR is effectively increased (or equivalently, the threshold $\sqrt{\rat}$ is lowered) in the PCA for the transformed matrix. The change of the threshold can be rigorously proved; see Theorem \ref{thm:trans-mean} for a precise statement. We remark a similar idea was also discussed in \cite{montanari2018adapting} without rigorous proof.

The corresponding result is not known, to our best knowledge, for the multiplicative model of the form $Y=(I+ \lambda \bsu \bsu^T)^{1/2} X =: (I+ \gamma \bsu \bsu^T) X$. (Here, $\lambda = 2\gamma + \gamma^2$.) The analysis is significantly more involved in this case due to the following reason: When applying a function $\f$ entrywise to $\sqrt{N} Y$, we find that
\[ \begin{split}
\f(\sqrt{N} Y_{ij} ) = \f \Big(\sqrt{N} X_{ij} + \gamma\sqrt{ N} \sum_k u_i u_k X_{kj} \Big) 
&\approx \f(\sqrt{N} X_{ij}) + \gamma \sqrt{N} \f'(\sqrt{N} X_{ij}) \sum_k u_i u_k X_{kj} \\
&\approx \sqrt{N} \left( \frac{\f(\sqrt{N} X_{ij})}{\sqrt{N}} + \gamma \E[\f'(\sqrt{N} X_{ij})] \sum_k u_i u_k X_{kj} \right),
\end{split} \]
and the transformed matrix is of the form $\gamma' \bsu \bsu^T X + \F$, which is not a spiked rectangular matrix anymore. (Note that $\F$ depends on $X$ and thus it cannot be considered as an additive model, either.)

In Theorem \ref{thm:trans-cov} in Section \ref{subsec:PCA}, we prove the effective change of the SNR for the multiplicative model. The proof of Theorem \ref{thm:trans-cov} is based on a generalized version of the BBP transition that works with the matrix of the form $\gamma \bsu \bsu^T X +\F$. Applying various results and techniques from random matrix theory, we introduce a general strategy to prove a BBP-type transition and apply it to the transformed matrix. 

It is notable that the optimal entrywise transform is different from the one for the additive model. For the additive model, the optimal transform is given by $-g'/g$, where $g$ is the density function of the noise entry. However, for the multiplicative model, the optimal transform is a linear combination of the function $-g'/g$ and the identity mapping. Heuristically, it is due to that the effective SNR depends not only on $\gamma'$ but also on the correlation between $X$ and $\F$; the former is maximized when the transform is $-g'/g$ while the latter is maximized when the transform is the identity mapping. We also remark that the effective SNR after the optimal entrywise transform is larger in the additive model, which suggests that the detection problem is fundamentally harder for the multiplicative model.

When it is impossible to reliably detect the signal, the next goal is the weak detection, which is basically the hypothesis testing problem between the null model and the alternative model that the spike exists in the data. As predicted by the Neyman--Pearson lemma, the likelihood ratio (LR) test is optimal in the sense that it minimizes the sum of the Type-I error and the Type-II error. The limit of the log-LR was proved to be Gaussian for both the additive model and the multiplicative model with Gaussian noise \cite{Onatski2013,el2018detection} from which the limiting optimal error can be readily deduced.

However, LR tests require substantial information on the prior, which is not available in many applications. Following the idea in \cite{chung2019weak}, we propose a test based on the LSS, which does not require any knowledge on the spike or the noise. We prove in Corollary \ref{cor:error} (see also Remark \ref{rem:optimal}) that the error of the proposed test is optimal if the noise is Gaussian. 

The proposed test is applicable even when the noise is non-Gaussian. It is expected that the weak detection based on the proposed test will perform better after the entrywise transform, which was proved for spiked Wigner models \cite{chung2019weak}. This will be discussed in a future paper. We also conjecture that with the entrywise transform our test will be optimal when the noise is non-Gaussian, but it is beyond our scope as the optimal error of the weak detection for non-Gaussian noise is not known, even for spiked Wigner models.

\subsection{Related works}

Spiked rectangular model was introduced by Johnstone \cite{Johnstone2001}. The transition of the largest eigenvalue was proved by Baik, Ben Arous, and P\'ech\'e \cite{BBP2005} for spiked complex Wishart matrices and generalized by Benaych-Georges and Nadakuditi \cite{Raj2011,benaych2012singular}. For more results from random matrix theory about the largest eigenvalue and the corresponding eigenvector of a spiked rectangular matrix, we refer to \cite{BKYY2016} and references therein.

The testing problem for spiked Wishart matrices with the spherical prior and Gaussian noise was considered by Onatski, Moreira, and Hallin \cite{Onatski2013,Onatski2014}, where they proved the optimal error of the hypothesis test. It is later extended to the case where the entries of the spikes are i.i.d. with bounded support (i.i.d. prior) by El Alaoui and Jordan \cite{el2018detection}.

The improved PCA based on the entrywise transformation was considered for spiked Wigner models in \cite{lesieur2015mmse,Perry2018}, where the transformation is chosen to maximize the effective SNR of the transformed matrix. Detection problems for spiked Wigner models were also considered, where the analysis is typically easier due to its symmetry and canonical connection with spin glass models. For more results on the spiked Wigner models, we refer to \cite{Montanari2017,Perry2018,AlaouiJordan2018,chung2019weak} and references therein.

\subsection{Organization of the paper}

The rest of the paper is organized as follows. In Section \ref{sec:prelim}, we precisely define the model and introduce previous results. In Section \ref{sec:main}, we state our results on the improved PCA and illustrate the improvement of PCA by numerical experiments. In Section \ref{sec:main2}, we state our results on the hypothesis testing and the central limit theorems for the linear spectral statistics. We conclude the paper in Section \ref{sec:conclusion} with the summary of our results and future research directions. Details of the numerical simulations and the proofs of the technical results can be found in Appendix.

\section{Preliminaries} \label{sec:prelim}

\subsection{Definition of the model}

We begin by precisely defining the model we consider in this paper. The noise matrix has the following properties.

\begin{defn}[Rectangular matrix] \label{defn:rect}
	We say an $M \times N$ random matrix $X = (X_{ij})$ is a (real) rectangular matrix if $X_{ij}$ ($1\leq i \leq M$, $1 \leq j\leq N$) are independent real random variables satisfying the following conditions:
	\begin{itemize}
		\item For all $i, j$, $\E[X_{ij}] = 0$, $N \E[X_{ij}^2]=1$, $N^{\frac{3}{2}} \E[X_{ij}^3]=w_3$, and $N^2 \E[X_{ij}^4]=w_4$ for some constants $w_3, w_4$. 
		\item For any positive integer $p$, there exists $C_p$, independent of $N$, such that $N^{\frac{p}{2}} \E[X_{ij}^p] \leq C_p$ for all $i, j$. 
	\end{itemize}
\end{defn}

The spiked rectangular matrices are defined as follows.

\begin{defn}[Spiked rectangular matrix - additive model] \label{defn:rect_mean}
	We say an $M \times N$ random matrix $Y = \sqrt{\lambda} \bsu \bsv^T + X$ is a rectangular matrix with spiked mean $\bsu$, $\bsv$ and SNR $\lambda$ if $\bsu = (u_1, u_2, \dots, u_M)^T \in \R^M$, $\bsv = (v_1, v_2, \dots, v_N)^T \in \R^N$ with $\| \bsu \| = \| \bsv \| = 1$, and $X$ is a rectangular matrix.
\end{defn}

\begin{defn}[Spiked rectangular matrix - multiplicative model] \label{defn:rect_cov}
	We say an $M \times N$ random matrix $Y = (I+ \lambda \bsu \bsu^T)^{1/2} X$ is a rectangular matrix with spiked covariance $\bsu$ and SNR $\lambda$ if $\bsu = (u_1, u_2, \dots, u_M)^T \in \R^M$ with $\| \bsu \|=1$ and $X$ is a rectangular matrix.
\end{defn}

We assume throughout the paper that $\lambda \geq 0$ and $\frac{M}{N} \to \rat \in (0, \infty)$ as $M, N \to \infty$.

\subsection{Principal component analysis} \label{subsec:BBP}

Let $S = YY^T$ be the sample covariance matrix (Gram matrix) derived from a spiked rectangular matrix $Y$. The empirical spectral measure of $S$ converges to the Marchenko--Pastur law $\mu_{MP}$, i.e., if we denote by $\mu_1 \geq \mu_2 \geq \dots \mu_M$ the eigenvalues of $S$, then
\beq \label{eq:MP_law}
\frac{1}{M} \sum_{i=1}^M \delta_{\mu_i}(x) \dd x \to \dd \mu_{MP}(x)
\eeq
weakly in probability as $M, N \to \infty$, where for $M \leq N$
\beq
\dd \mu_{MP}(x) = \frac{\sqrt{(x- d_-)(d_+ -x)}}{2\pi\rat x} \mathbf{1}_{(d_-, d_+)}(x) \dd x,
\eeq
with $d_{\pm} = (1\pm \sqrt{\rat})^2$.
The largest eigenvalue has the following (almost sure) limit:
\begin{itemize}
	\item If $\lambda > \sqrt{\rat}$, then $\mu_1 \to (1+\lambda)(1+\frac{\rat}{\lambda})$.
	\item If $\lambda < \sqrt{\rat}$, then $\mu_1 \to d_+ = (1+\sqrt{\rat})^2$.
\end{itemize}

This in particular shows that the detection can be reliably done by PCA if $\lambda > \sqrt{\rat}$.

\subsection{Likelihood ratio} \label{subsec:LR}

Denote by $\p_1$ the joint probability of the data $Y$, a spiked rectangular matrix, with $\lambda =\SNR> 0$ and $\p_0$ with $\lambda = 0$. When the noise is Gaussian, the likelihood ratio $\caL(Y; \lambda)$ of $\p_1$ with respect to $\p_0$ is given by
\[ \begin{split}
&\int \frac1{\det(I+\SNR \bsu \bsu^T)^{\frac{N}{2}}}\cdot\exp\left(\frac{N\lambda}{2(1+\SNR \| \bsu \|^2)} \sum_{i=1}^{M}\sum_{j=1}^{N} (YY^T)_{ij} u_i u_j\right) \dd \caP_{\bsu}
\end{split} \]
for the multiplicative model (Definition \ref{defn:rect_cov}) and
\[ \begin{split}
\int\exp\left(N\sum_{i=1}^{M}\sum_{j=1}^{N} \left[ \sqrt{\SNR} Y_{ij} u_i v_j - \frac{\SNR}{2} u_i^2 v_j^2 \right] \right) \dd \caP_{\bsu} \dd \caP_{\bsv}
\end{split} \]
for the additive model (Definition \ref{defn:rect_mean}). Here, $\caP_{\bsu}$ and $\caP_{\bsv}$ are the prior distributions of $\bsu$ and $\bsv$, respectively. 

If $\SNR < \sqrt{\rat}$, for both models with the spherical prior where the spike is drawn uniformly from the unit sphere, the log-LR has the Gaussian limit; as $N \to \infty$, it converges to
\[ \begin{split}
\caN \left( \frac{1}{4} \log\left(1-\frac{\SNR^2}{\rat}\right), -\frac{1}{2}\log\left(1-\frac{\SNR^2}{\rat}\right)\right)
\end{split} \]
under the null hypothesis $\bsH_0: Y \sim \p_0$ and
\[ \begin{split}
\caN \left( -\frac{1}{4} \log\left(1-\frac{\SNR^2}{\rat}\right), -\frac{1}{2}\log\left(1-\frac{\SNR^2}{\rat}\right)\right)
\end{split} \]
under the alternative hypothesis $\bsH_1: Y \sim \p_1$. The same result also holds for the additive model with Rademacher prior. The sum of the Type-I error and the Type-II error of the likelihood ratio test
\begin{align} 
\label{eq:LR_error}
\err(\SNR) &:= \p(L(Y;\SNR) > 1 | \bsH_0) + \p(L(Y;\SNR) \leq 1 | \bsH_1) \to \erfc \left( \frac{1}{4} \sqrt{-\log \left( 1- \frac{\SNR^2}{\rat} \right)} \right)
\end{align}
as $N \to \infty$. We remark that it is the minimal error among all tests as Neyman--Pearson lemma asserts. This in particular shows that the reliable detection of signal is impossible with Gaussian noise when $\SNR < \sqrt{\rat}$.

\subsection{Linear spectral statistics}

The proof of the Gaussian convergence of the LR in \cite{Baik-Lee2016,Baik-Lee2018} is based on the recent study of linear spectral statistics, defined as
\beq \label{eq:LSS}
L_Y(f) = \sum_{i=1}^M f(\mu_i)
\eeq
for a function $f$, where $\mu_1 \geq \mu_2 \geq \dots \mu_M$ are the eigenvalues of $S=YY^T$. As the Marchenko--Pastur law in \eqref{eq:MP_law} suggests, it is required to consider the fluctuation of the LSS about
\[
M \int_{d_-}^{d_+} f(x) \, \dd \mu_{MP}(x).
\]
The CLT for the LSS is the statement
\beq \begin{split} \label{eq:CLT_LSS}
	&\left( L_Y(f) - M \int_{d_-}^{d_+} f(x) \, \dd \mu_{MP}(x) \right) \Rightarrow \caN(m_Y(f), V_Y(f)),
\end{split} \eeq
where the right-hand side is the Gaussian random variable with the mean $m_Y(f)$ and the variance $V_Y(f)$. The CLT was proved for the null case ($\lambda=0$). We will show that the CLT also holds under the alternative and the mean $m_Y(f)$ depends on $\lambda$ while the variance $V_Y(f)$ does not.

\section{Main Result I - Improved PCA} \label{sec:main}

In this section, we state our first main results on the improvement of PCA by entrywise transformations and provide the results from numerical experiments.

\subsection{Improved PCA} \label{subsec:PCA}

Let $\caP$ be the distribution of the normalized entry $\sqrt{N} X_{ij}$ whose density function is $g$. 
As we discussed in Section \ref{subsec:main}, applying a function $\f$ to the additive model in Definition \ref{defn:rect_cov} approximately yields another rectangular matrix
\beq \label{eq:approx_mean}
\frac{\f(\sqrt{N} X_{ij})}{\sqrt{N}} + \sqrt{\lambda} \E[\f'(\sqrt{N} X_{ij})] u_i v_j.
\eeq
Suppose that $\f(\sqrt{N} X_{ij}) = \sqrt{N} \F_{ij}$ is with mean $0$ and variance $1$. Then, the effective SNR of the transformed matrix is $\lambda (\E[\f'(\sqrt{N} X_{ij})])^2$, which is maximized when $\f(x)$ is a multiple of $-g'(x)/g(x)$.

For the multiplicative model in Definition \ref{defn:rect_cov}, applying a function $\f$ approximately yields a transformed matrix of the form $\F + \wh \gamma \bsu \bsu^T X$ as discussed in Section \ref{subsec:main}, where we set $\wh \gamma = \gamma \E[\f'(\sqrt{N} X_{ij})]$. The sample covariance matrix generated by it is
\[ \begin{split}
&(\F + \wh \gamma \bsu \bsu^T X)(\F + \wh \gamma \bsu \bsu^T X)^T = \F\F^T + \wh \gamma \F X^T \bsu \bsu^T + \wh \gamma \bsu \bsu^T X \F^T + \wh \gamma^2 \bsu \bsu^T X X^T \bsu \bsu^T.
\end{split} \]
Conditioning on $\bsu$, its expectation is $(I + \lambda_{SNR} \bsu \bsu^T)$, where the effective SNR $\lambda_{SNR}$ is
\[ \begin{split}
&2\wh \gamma \E[\sqrt{N} X_{ij} \f(\sqrt{N} X_{ij})] + \wh \gamma^2 =2\gamma \E[\f'(\sqrt{N} X_{ij})] \E[\sqrt{N} X_{ij} \f(\sqrt{N} X_{ij})] + \gamma^2 (\E[\f'(\sqrt{N} X_{ij})])^2
\end{split} \]
We can find that $\lambda_{SNR}$ is maximized when $\f(x)$ is a multiple of $-g'(x)/g(x) + \alpha x$ for some constant $\alpha$.

In this section, we rigorously prove our heuristic argument and show the detection threshold of PCA can be lowered by applying the entrywise transformations above.
We introduce the following assumptions for the spike and the noise.
\begin{assump} \label{assump:entry}
	For the spike $\bsu$ (and also $\bsv$ in the additive model), we assume either 
	\begin{enumerate}
		\item the spherical prior, i.e., $\bsu$ (and $\bsv$) are drawn uniformly from the unit sphere, or
		\item the i.i.d. prior, i.e., the entries $u_1, u_2, \dots, u_M$ (respectively, $v_1, v_2, \dots, v_N$) are i.i.d. random variables with mean zero and variance $M^{-1}$ (respectively $N^{-1}$) such that for any integer $p>2$
		\[
		\E |u_i|^p, \E |v_j|^p \leq \frac{C_p}{M^{1+(p-2)\phi}}
		\]
		for some ($N$-independent) constants $C_p>0$ and $\phi > \frac{1}{4}$, uniformly on $i$ and $j$.
	\end{enumerate} 	
	
	For the noise, let $\caP$ be the distribution of the normalized entries $\sqrt{N} X_{ij}$. We assume the following:
	\begin{enumerate}
		\item The density function $g$ of $\caP$ is smooth, positive everywhere, and symmetric (about 0).
		\item For any fixed $D$, the $D$-th moment of $\caP$ is finite.
		\item The function $h = -g'/g$ and its all derivatives are polynomially bounded in the sense that $|h^{(\ell)}(w)| \leq C_{\ell} |w|^{C_{\ell}}$ for some constant $C_{\ell}$ depending only on $\ell$.
	\end{enumerate}  
\end{assump}

Note that the signal is not necessarily delocalized, i.e., some entries of the signal can be much larger than $N^{-1/2}$.

We remark that some conditions in Assumption \ref{assump:entry}, especially the i.i.d. prior and the finiteness of all moments of $\caP$, are technical constraints and our results hold under weaker assumptions. We also remark that if $\sqrt{M} u_i$ (and $\sqrt{N} v_j$) are i.i.d. random variables, independent of $M$ (and $N$), whose all moments are finite, Assumption \ref{assump:entry} is satisfied with $\phi = \frac{1}{2}$.

Given the data matrix $Y$, we consider a family of the entrywise transformations of the form $h_{\alpha}(x) = -g'(x)/g(x) + \alpha x$ and transformed matrices $\wt Y^{(\alpha)}$ whose entries are
\beq \label{eq:transformed-cov}
\wt Y^{(\alpha)}_{ij} = \frac{1}{\sqrt{(\alpha^2+2\alpha+\fg) N}} h_{\alpha}(\sqrt{N} Y_{ij}),
\eeq
where the Fisher information $\fg$ of $g$ is given by
\[
\fg = \int_{-\infty}^{\infty} \frac{(g'(x))^2}{g(x)} \dd x.
\]
Note that $\fg \geq 1$ where the equality holds only if $g$ is the standard Gaussian.

For the additive model, we show that the effective SNR of the transformed matrix for PCA is $\lambda \fg$.
\begin{thm} \label{thm:trans-mean}
	Let $Y$ be a spiked rectangular matrix in Definition \ref{defn:rect_mean} that satisfy Assumption \ref{assump:entry}. Let $\wt Y \equiv \wt Y^{(0)}$ be the transformed matrix obtained as in \eqref{eq:transformed-cov} with $\alpha=0$ and $\wt \mu_1$ the largest eigenvalue of $\wt Y \wt Y^T$. Then, almost surely,
	\begin{itemize}
		\item If $\lambda > \frac{\sqrt{\rat}}{\fg}$, then $\wt\mu_1 \to (1+\lambda \fg)(1+\frac{\rat}{\lambda \fg})$.
		\item If $\lambda < \frac{\sqrt{\rat}}{\fg}$, then $\wt\mu_1 \to d_+ = (1+\sqrt{\rat})^2$.
	\end{itemize}
\end{thm}

From Theorem \ref{thm:trans-mean}, if $\lambda > \frac{\sqrt{\rat}}{\fg}$, we immediately see that the signal in the additive model can be reliably detected by the transformed PCA. Thus, the detection threshold in the PCA is lowered when the noise is non-Gaussian. We also remark that $h_0$ is the optimal entrywise transformation (up to constant factor) as in the Wigner case; see Appendix \ref{sec:optimize_transform}.

For the proof, we first adapt the strategy in \cite{Perry2018} to justify that the transformed matrix is approximately equal to \eqref{eq:approx_mean}, which is another rectangular matrix. We then prove a BBP-type transition for the additive model, following the method of \cite{benaych2012singular}. Since our assumptions on the spike and the noise are weaker, we provide the detail of the proof of Theorem \ref{thm:trans-mean} in Appendix \ref{sec:proof_trans-mean}.

For the multiplicative model, we have the following.

\begin{thm} \label{thm:trans-cov}
	Let $Y$ be a spiked rectangular matrix in Definition \ref{defn:rect_cov} that satisfy Assumption \ref{assump:entry}. Let $\wt Y \equiv \wt Y^{(\alpha_g)}$ be the transformed matrix obtained as in \eqref{eq:transformed-cov} with 
	\[
	\alpha_g := \frac{-\gamma F_g+\sqrt{4F_g+4\gamma F_g+\gamma^2F_g^2}}{2(1+\gamma)}
	\]
	and $\wt \mu_1$ the largest eigenvalue of $\wt Y \wt Y^T$. Then, almost surely,
	\begin{itemize}
		\item If $\lambda_g > \sqrt{\rat}$, then $\wt\mu_1 \to (1+ \lambda_g)(1+\frac{\rat}{\lambda_g})$.
		\item If $\lambda_g < \sqrt{\rat}$, then $\wt\mu_1 \to d_+ = (1+\sqrt{\rat})^2$.
	\end{itemize}
	where
	\[
	\lambda_g := \gamma + \frac{\gamma^2 \fg}{2} + \frac{\gamma \sqrt{4\fg + 4\gamma \fg + \gamma^2 \fg^2}}{2}.
	\]
\end{thm}

Note that
\[ \begin{split}
\lambda_g &\geq \gamma + \frac{\gamma^2 \fg}{2} + \frac{\gamma \sqrt{4 + 4\gamma \fg + \gamma^2 \fg^2}}{2} = 2\gamma + \gamma^2 F_g \geq 2\gamma + \gamma^2 = \lambda,
\end{split} \]
and the inequality is strict if $F_g > 1$, i.e., $g$ is not Gaussian. 
From Theorem \ref{thm:trans-cov}, if $\lambda_g > \sqrt{\rat}$, we immediately see that the signal can be reliably detected by the transformed PCA. Thus, the detection threshold in the PCA is lowered when the noise is non-Gaussian. We also remark that $h_{\alpha_g}$ is the optimal entrywise transformation (up to constant factor); see Appendix \ref{sec:optimize_transform}.

We outline the proof of Theorem \ref{thm:trans-cov}. We begin by justifying that the transformed matrix $\wt Y$ is approximately of the form $(\F + \wh\gamma \bsu \bsu^T X)$. Then, the largest eigenvalue of $\wt Y \wt Y^T$ can be approximated by the largest eigenvalue of $(\F + \wh\gamma \bsu \bsu^T X)^T (\F + \wh\gamma \bsu \bsu^T X)$ for which we consider an identity
\[ \begin{split}
&(\F + \wh\gamma \bsu \bsu^T X)^T (\F + \wh\gamma \bsu \bsu^T X) - zI = (\F^T \F -zI) (I + L(z)),
\end{split} \]
where
\begin{align*}
&L(z) = \caG(z) (\wh \gamma X^T \bsu \bsu^T \F + \wh \gamma \F^T \bsu \bsu^T X + \wh \gamma^2 X^T \bsu \bsu^T X), &&\caG(z) = (\F^T \F -zI)^{-1}.
\end{align*}
If $z$ is an eigenvalue of $(\F + \wh\gamma \bsu \bsu^T X)^T (\F + \wh\gamma \bsu \bsu^T X)$ but not of $\F^T \F$, the determinant of $(I + L(z))$ must be $0$ and hence $-1$ is an eigenvalue of $L(z)$. Since the rank of $L(z)$ is at most $2$, we can find that the eigenvector of $L(z)$ is a linear combination of two vectors $\caG(z) \F^T \bsu$ and $\caG(z) X^T \bsu$, i.e., for some $a, b$,
\beq \begin{split} \label{eq:eigenvalue}
	&L(z) (a \caG(z) \F^T \bsu + b \caG(z) X^T \bsu) = -(a \caG(z) \F^T \bsu + b \caG(z) X^T \bsu).
\end{split} \eeq

From the definition of $L(z)$,
\[ \begin{split}
L(z) \cdot \caG(z) X^T \bsu &= \wh \gamma \langle \bsu, \F \caG(z) X^T \bsu \rangle \cdot \caG(z) X^T \bsu  + \wh \gamma \langle \bsu,X \caG(z) X^T \bsu \rangle \cdot \caG(z) \F^T \bsu \\& \quad + \wh \gamma^2 \langle \bsu,X \caG(z) X^T \bsu \rangle \cdot \caG(z) X^T \bsu,
\end{split} \]
and a similar equation holds for $L(z) \cdot \caG(z) \F^T \bsu$. It suggests that if $\langle \bsu, \F \caG(z) X^T \bsu \rangle$ and $\langle \bsu,X \caG(z) X^T \bsu \rangle$ are concentrated around deterministic functions of $z$, then the left side of \eqref{eq:eigenvalue} can be well-approximated by a (deterministic) linear combination of $\caG(z) \F^T \bsu$ and $\caG(z) X^T \bsu$. We can then find the location of the largest eigenvalue in terms of a deterministic function of $z$ and conclude the proof by optimizing the function $\f$.

The concentration of random quantities $\langle \bsu, \F \caG(z) X^T \bsu \rangle$ and $\langle \bsu,X \caG(z) X^T \bsu \rangle$ is the biggest technical challenge in the proof, mainly due to the dependence between the matrices $\F$ and $X$. We prove it by applying the technique of linearization in conjunction with resolvent identities and also several recent results from random matrix theory, most notably the local Marchenko--Pastur law. 

The detailed proof of Theorem \ref{thm:trans-cov} can be found in Appendix \ref{sec:proof_trans-cov}.

\begin{rem}
	Unlike the additive model, we cannot determine $\alpha_g$ without prior knowledge on the SNR. Nevertheless, we can apply the transformation $h_{\sqrt{\fg}}$, which effectively increases the SNR; see Appendix \ref{sec:optimize_transform}.
\end{rem}

\subsection{Applying the improved PCA to real data}

\begin{figure}[t]
	\vskip 0.2in
	\begin{center}
		\centerline{\includegraphics[width=0.7\columnwidth]{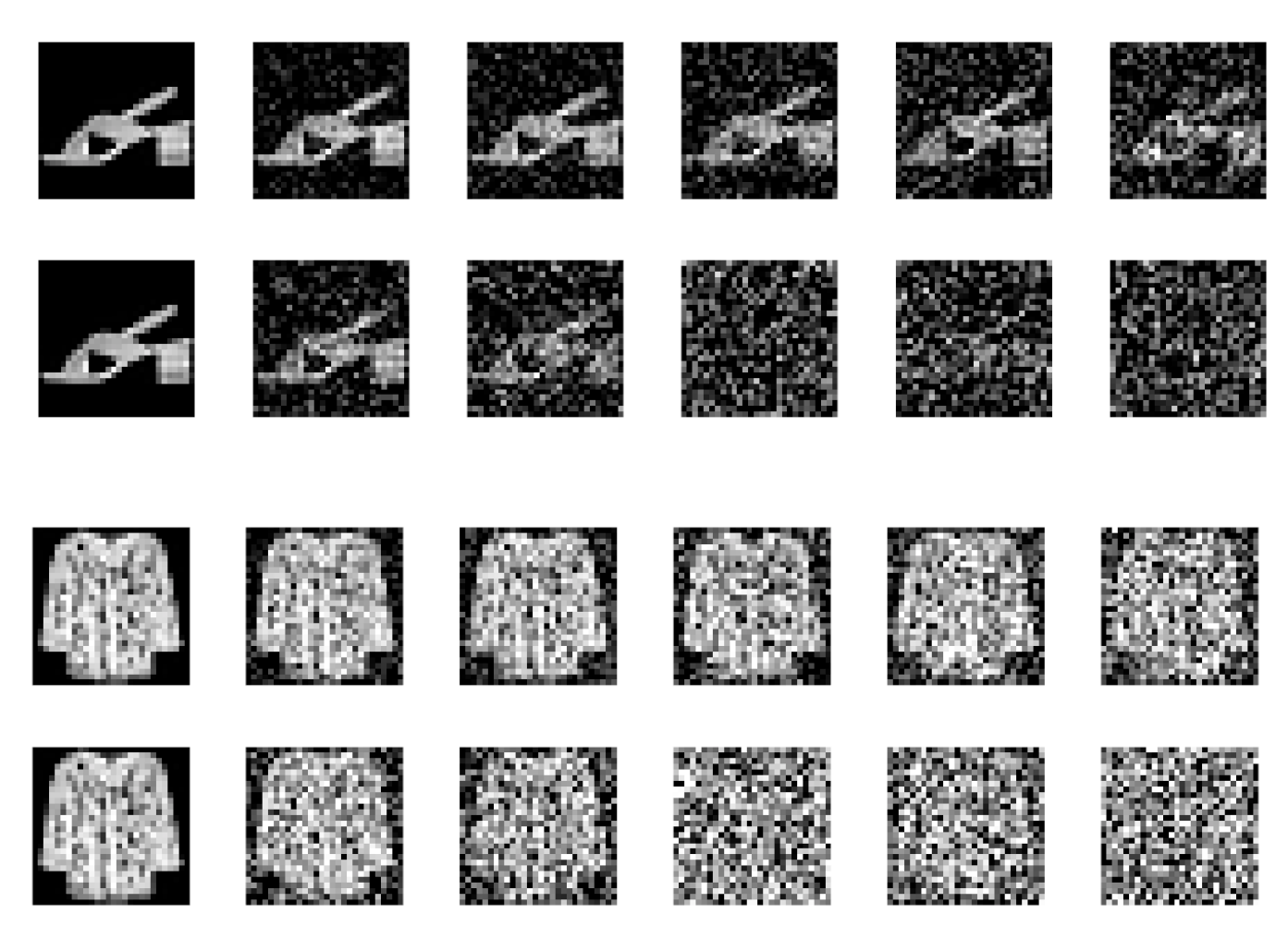}}
		\caption{We compare the reconstruction performance of the proposed PCA (top lines) and the standard PCA (bottom lines) for two Fashion-MNIST images, with the number of measurements $N=[3136,1568,784,588,392]$ and $M=784$. The left most column displays the original images for comparison.}
		\label{fig:real}
	\end{center}
	\vskip -0.2in
\end{figure}

To illustrate the improvement of PCA in Section \ref{subsec:PCA}, we perform the following numerical experiment:
We choose a vector $\bsz \in \R^{784}$ from the standard Fashion-MNIST dataset. We then let the spike $\bsu$ be a normalized vector of $\bsz$.
The $j$-th column of the data matrix $Y$ is a noisy sample of the spike given by
\[
Y_j = v_j \bsu + X_j,
\]
where $v_j$ follows Rademacher distribution and each entry of $X_j$ is independently drawn from a centered bimodal distribution with unit variance, which is a convolution of Gaussian and Rademacher random variables, and normalized by $1/\sqrt{N}$. 
Our goal is to reconstruct the spike $\bsu$ from $Y$ with $N$ columns. In Fig.~\ref{fig:real}, we compare the reconstruction by the improved PCA with standard PCA over $Y$. With the optimal entrywise transformation, the proposed PCA outperforms the standard PCA.

While we have analyzed the improved PCA with prior information on the noise, it is possible to estimate the noise even when the noise distribution is not known. As an attempt, we tried kernel density estimation (KDE) with the Gaussian kernel, which approximates the density of the noise $g(x)$ by
\[
	\widehat{g}(x) := \frac{1}{MN\delta} \sum_{i, j} \phi( (x-\sqrt{N}Y_{ij}) / \delta),
\]
where $\phi$ is the density function of the standard normal random variable and $\delta$ is the bandwidth, which we chose to be $(MN)^{-1/5}$. 

For a numerical experiment, we consider the data matrix $Y= \sqrt{\lambda} \bsu \bsv^T + X$, where $\sqrt{M} u_i$ and $\sqrt{N} v_j$ follow Rademacher distribution for $i=1, 2, \dots, M$ and $j=1, 2, \dots, N$. The noise is independently drawn from the same centered bimodal distribution as in the experiment above but with the variance $N^{-1}$. The size of the data matrix is set to be $M=1024$, $N=2048$, and hence the ratio $\rat = M/N = 1/2$. We set the SNR $\lambda \approx 0.4945$. With the approximation $\widehat{g}$, we use the entrywise transformation $\widehat{h}:= -\widehat{g}'/\widehat{g}$.

\begin{figure}[t]
	\vskip 0.2in
	\begin{center}
		\centerline{\includegraphics[width=\columnwidth]{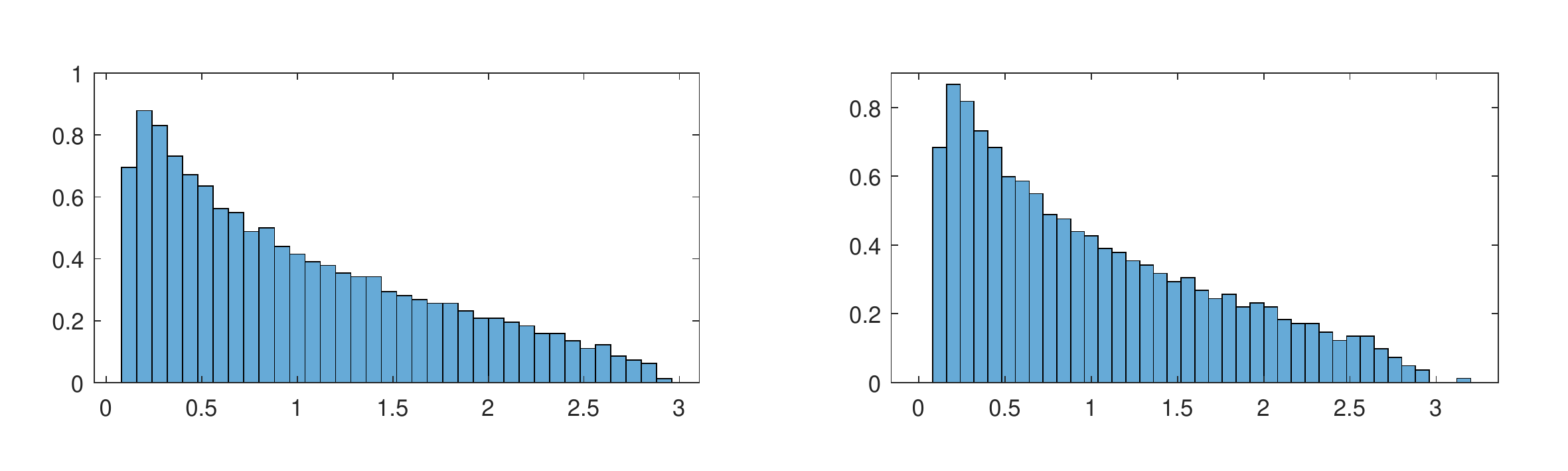}}
		\caption{The spectrum of the sample covariance matrices, before (left) and after (right) the entrywise transformation $\widehat{h}$. An outlier eigenvalue can be seen only after the entrywise transformation}
		\label{fig:KDE}
	\end{center}
	\vskip -0.2in
\end{figure}

In Fig.~\ref{fig:KDE}, we compare the spectrum of the sample covariance matrices, $YY^T$ (left) and $\wt Y \wt Y^T$ (right), where for the latter we rescale the eigenvalues so that the bulk of its spectrum matches that of the former.
An isolated eigenvalue can be seen only in the spectrum in the bottom, which is the case after the entrywise transformation.

For more simulation results about the improved PCA, see Appendix \ref{sec:simulation}. 

\section{Main Result II - Weak Detection} \label{sec:main2}

In this section, we state our second main results on the hypothesis test and provide the results from numerical experiments.

\subsection{Hypothesis testing and central limit theorem} \label{subsec:ht}

Suppose that the SNR $\SNR$ for the alternative hypothesis $\bsH_1$ is known and our goal is to detect the presence of the signal. We propose a test based on the LSS of the data matrix in \eqref{eq:LSS}. 
The key observation is that the variances of the limiting Gaussian distributions of the LSS are equal while the means are not. If we denote by $V_Y(f)$ the common variance, and $m_Y(f)|_{\bsH_0}$ and $m_Y(f)|_{\bsH_1}$ the means, respectively, our goal is to find a function that maximizes the relative difference between the limiting distributions of the LSS under $\bsH_0$ and under $\bsH_1$, i.e.,
\[
\left| \frac{m_Y(f)|_{\bsH_1} - m_Y(f)|_{\bsH_0}}{\sqrt{V_Y(f)}} \right|.
\]
As we will see in Theorem \ref{thm:CLT}, the optimal function $f$ is of the form $C_1\phi_{\SNR} +C_2$ for some constants $C_1$ and $C_2$, where 
\beq \begin{split}\label{eq:optimal_f}
	\phi_{\SNR}(x) &=\frac{\SNR}{\rat}\left(\frac{2}{w_4-1}-1\right)x  -\log\left(\left(1+\frac{\rat}{\SNR} \right)(1+\SNR) - x\right).
\end{split}\eeq
The test statistic we use is thus defined as
\beq\label{eq:L_lambda}
\begin{split}
	L_{\SNR} &= \sum_{i=1}^M \phi_{\SNR}(\mu_i) - M \int_{d_-}^{d_+} \phi_{\SNR}(x) \,\dd \mu_{MP}(x) \\
	&= -\log \det \left( \left(1+\frac{\rat}{\SNR} \right)(1+\SNR)I - YY^T \right)  + \frac{\SNR}{\rat} \left( \frac{2}{w_4-1} - 1 \right) (\Tr YY^T-M) \\
	&\quad + M \left[\frac{\SNR}{\rat} - \log\left(\frac{\SNR}{\rat}\right) -\frac{1-\rat}{\rat}\log(1+\SNR) \right].
\end{split} \eeq

Our main result in this section is the following CLT for $L_{\SNR}$.
\begin{thm} \label{thm:test}
	Let $Y$ be a spiked rectangular matrix in Definition \ref{defn:rect_cov} or \ref{defn:rect_mean} with $w\in(0,\sqrt{\rat})$ and $w_4>1$. Then, for any spikes with $\|\bsu\|=1$ and $\|\bsv\|=1$,
	\beq
	L_{\SNR} \Rightarrow \mathcal{N}(m(\lambda),V_0).
	\eeq
	The mean of the limiting Gaussian distribution is given by
	\beq \begin{split} \label{eq:mean_test}
		m(\lambda) &= -\frac{1}{2} \log\left(1-\frac{\SNR^2}{\rat}\right)  -\frac{\SNR^2}{2\rat}(w_4-3)  -\log\left(1-\frac{\lambda^2}{\rat}\right) +\frac{\lambda^2}{\rat}\left(\frac{2}{w_4-1}-1\right)
	\end{split} \eeq
	and the variance
	\beq \label{eq:var_test}
	V_0=-2\log\left(1-\frac{\SNR^2}{\rat}\right) +\frac{2\SNR^2}{\rat}\left(\frac2{w_4-1}-1\right)
	\eeq
\end{thm}

Theorem \ref{thm:test} directly follows from the general CLT result in Theorem \ref{thm:CLT}. See also Appendix \ref{sec:compute}
for more detail on the mean and the variance.

We propose a test in Algorithm \ref{alg:test} based on Theorem \ref{thm:test}. In this test, we compute the test statistic $L_{\SNR}$ and compare it with the average of $m(0)$ and $m(\SNR)$, i.e.,
\beq \begin{split} \label{eq:m_lambda}
	m_{\SNR} &:= \frac{m(0)+m(\SNR)}{2} = -\log\left(1-\frac{\SNR^2}{\rat}\right) +\frac{\SNR^2}{2\rat} \left( \frac{2}{w_4-1} -w_4 +2 \right).
\end{split} \eeq

\begin{algorithm}[tb]
	\caption{Hypothesis test}
	\label{alg:test}
	\begin{algorithmic}
		\STATE {\bfseries Input:} data $Y_{ij}$, parameters $w_4$, $\SNR$
		\STATE $L_{\SNR} \gets$ test statistic in \eqref{eq:L_lambda}
		\STATE $m_{\SNR} \gets$ critical value in \eqref{eq:m_lambda}
		\IF{$L_{\SNR} \leq m_{\SNR}$}
		\STATE Accept $\bsH_0$
		\ELSE
		\STATE Reject $\bsH_0$
		\ENDIF
	\end{algorithmic}
\end{algorithm}

As a simple corollary to Theorem \ref{thm:test}, we have the following formula for the limiting error of the proposed test.

\begin{cor} \label{cor:error}
	The error of the test in algorithm \ref{alg:test},
	\beq \begin{split} \label{eq:test_error}
		\err(\SNR) &= \p( L_{\SNR} > m_{\SNR} | \bsH_0) + \p( L_{\SNR} \leq m_{\SNR} | \bsH_1) \to \erfc \left( \frac{\sqrt{V_0}}{4\sqrt 2} \right),
	\end{split} \eeq
	where $V_0$ is the variance in \eqref{eq:var_test} and $\erfc(\cdot)$ is the complementary error function.
\end{cor}

For the proof of Corollary \ref{cor:error}, see Corollary 5 of \cite{AlaouiJordan2018} or Theorem 2 of \cite{chung2019weak}.

\begin{rem} \label{rem:optimal}
If the noise $X$ is Gaussian, $w_4=3$ and the limiting error in Corollary \ref{cor:error} is
\[
\erfc \left( \frac{\sqrt{V_0}}{4\sqrt 2} \right) = \erfc \left( \frac{1}{4} \sqrt{-\log \left( 1- \frac{\SNR^2}{\rat} \right)} \right),
\]
and it coincides with the error of the LR test in \eqref{eq:LR_error}. It shows that our test is optimal with the Gaussian noise.
\end{rem}

Even if the exact parameter $w_4$ is not known a priori, it can be easily estimated from the data matrix $Y$ by computing $\frac{1}{MN} \sum Y_{ij}^4$. The accuracy of such an estimate can be easily checked from the Chernoff bound.

Lastly, we state a general CLT for the LSS and the optimality of the function $\phi_{\SNR}$ as the test statistic.
\begin{thm} \label{thm:CLT}
	Assume the conditions in Theorem \ref{thm:test}. Denote by $\mu_1 \geq \mu_2 \geq \cdots\geq\mu_M$ the eigenvalues of $YY^T$. For any function $f$ analytic on an open set containing an interval $[d_-,d_+],$
	\beq \begin{split} \label{eq:CLT_statement}
		&\left(\sum_{i=1}^{M}f(\mu_i)-M \int_{d_-}^{d_+} \frac{\sqrt{(x-d_-)(d_+ -x)}}{2\pi \rat x} \phi_{\SNR}(x) \, \dd x \right) \Rightarrow \mathcal{N}(m_Y(f),V_Y(f)).
	\end{split} \eeq
	The mean and the variance of the limiting Gaussian distribution are given by
	\[ \begin{split}
	m_Y(f) &= \frac{\wt{f}(2) + \wt{f}(-2)}{4} -\frac{\tau_0(\wt{f})}{2} -(w_4-3)\tau_2(\wt{f})  +\sum_{\ell=1}^{\infty}\left(\frac{\SNR}{\sqrt{\rat}}\right)^{\ell}\tau_\ell(\wt{f})
	\end{split} \]
	and
	\[
	V_Y(f) =2\sum_{\ell=1}^\infty \ell\tau_\ell(\wt{f})^2+(w_4-3)\tau_1(\wt{f})^2,
	\]
	where we let $\wt f(x)=f(\sqrt{\rat}x+1+\rat)$, 
	\[
	\tau_{\ell}(f) = \frac{1}{\pi} \int_{-2}^2 T_{\ell}\left( \frac{x}{2} \right) \frac{f(x)}{\sqrt{4-x^2}} \dd x\,,
	\]
	and $T_{\ell}$ is the $\ell$-th Chebyshev polynomial of the first kind.
	
	Furthermore, for $m(\SNR)$, $m(0)$, and $V_0$ defined in Theorem \ref{thm:test},
	\[
	\left|\frac{m_Y(f) - m_X(f)}{\sqrt{V_{Y}(f)}}\right| \leq \left|\frac{m(\SNR)- m(0)}{\sqrt{V_0}}\right|
	\]
	The equality holds if and only if $f(x)=C_1\phi_\SNR(x) +C_2$ for some constants $C_1$ and $C_2$ with the function $\phi_\SNR$ defined in \eqref{eq:optimal_f}.
\end{thm}

We remark that the analyticity of the function $f$ in Theorem \ref{thm:CLT} is assumed only because it is sufficient in our purpose and this assumption can be weakened by the density argument, which is typically used in the proof of CLT results in random matrix theory.

We briefly sketch the proof of Theorem \ref{thm:CLT} based on the interpolation technique, developed in \cite{chung2019weak,jung2020weak}. In this method, the right side of \eqref{eq:CLT_statement} is written as the following contour integral of the trace of the resolvent: For a function $f$ analytic on an open set containing an interval $[d_-, d_+]$,
\beq \label{eq:Cauchy} \begin{split}
	\sum_{i=1}^M f(\mu_i) &= \sum_{i=1}^M \frac{1}{2\pi\ii} \oint_{\Gamma} \frac{f(z)}{z-\mu_i} \dd z \\
	&= -\frac{1}{2\pi\ii} \oint_{\Gamma} f(z) \Tr (YY^T -zI)^{-1} \dd z
\end{split} \eeq
for any contour $\Gamma$ containing $\mu_1, \mu_2, \dots, \mu_N$. For the null model, i.e., if $\lambda = 0$, the CLT was proved in \cite{Bai-Silverstein2004,Lytova-Pastur2009a} with precise formulas for the mean and the variance. 

To prove the CLT for a non-null model, i.e., a spiked rectangular matrix with $\lambda \neq 0$, we introduce an interpolation between the null model and the non-null model, and track the change of the LSS by finding the change of $\Tr (YY^T -zI)^{-1}$. The change is decomposed into the deterministic part and the random part, where the latter converges to $0$ with overwhelming probability for both the additive model and the multiplicative model. We can then conclude that the CLT for the LSS holds also for the non-null model, with the variance unchanged. The change of the mean can be computed by considering the deterministic change of the resolvent. 

The detail of the proof of Theorem \ref{thm:CLT} can be found in Appendix \ref{sec:proof_CLT}.

\subsection{Numerical experiments for the LSS test} \label{subsec:ht_sim}

We consider the case where the noise matrix $X$ is Gaussian and the signal $\bsu = (u_1, u_2, \dots, u_M)^T$ and $\bsv = (v_1, v_2, \dots, v_N)^T$, where $\sqrt{M} u_i$'s and $\sqrt{N} v_j$'s are i.i.d. Rademacher random variables for $i=1, 2, \dots, M$ and $j=1, 2, \dots, N$. Let the data matrix $Y = \sqrt{\lambda} \bsu \bsv^T + X$. The parameters are $w_2 = 2$ and $w_4 = 3$.

\begin{figure}[t]
	\vskip 0.2in
	\begin{center}
		\centerline{\includegraphics[width=200pt]{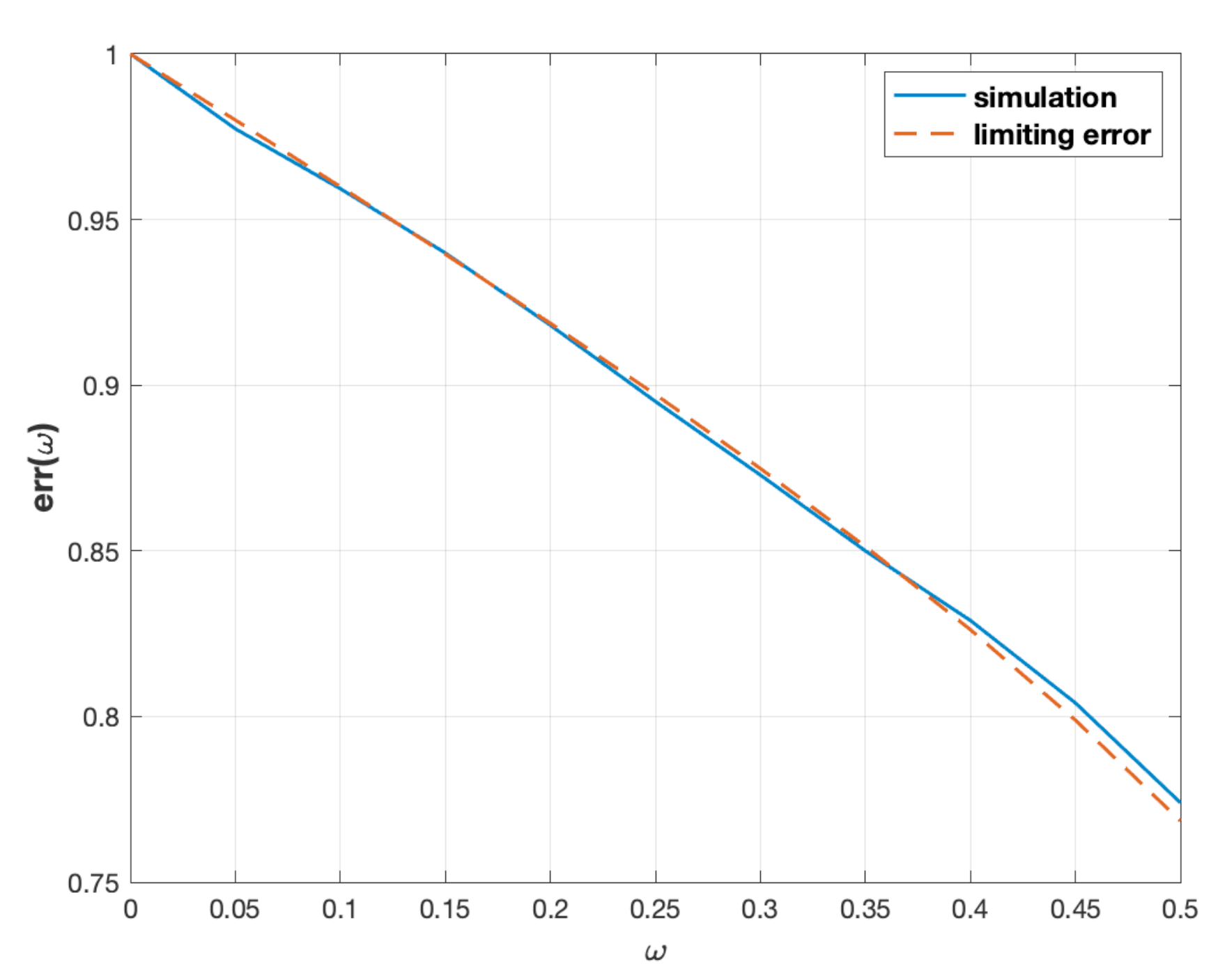}}
		\caption{The error from the simulation (solid) and the theoretical limiting error in \eqref{eq:test_error} (dashed), respectively, for the Gaussian noise.}
		\label{fig:GOE}
	\end{center}
	\vskip -0.2in
\end{figure}

In Figure \ref{fig:GOE}, we plot empirical average (after 10,000 Monte Carlo simulations) of the error of the proposed test in Algorithm \ref{alg:test} and the theoretical (limiting) error in \eqref{eq:test_error}, varying the SNR $\SNR$ from $0$ to $0.5$, with $M=256$ and $N=512$. It can be checked that the error of the proposed test closely matches the theoretical error.

\section{Conclusion and Future Works} \label{sec:conclusion}

In this paper, we considered the detection problem of spiked rectangular model. For both the multiplicative model and the additive model, we showed that PCA can be improved for non-Gaussian noise by transforming the data entrywise. We proved the effective SNR and the optimal entrywise transforms for both models. We also proposed a universal test that does not require any prior information on the spike. The test and its error do not depend on the noise except its (normalized) fourth moment. The error of the proposed test is optimal when the noise is Gaussian.

A natural future research direction is to apply the entrywise transformation for the weak detection. As in the spiked Wigner model, we believe that the error of the proposed test can be lowered with the entrywise transformation and it can be proved by establishing the central limit theorems for the transformed matrices. 

\section*{Acknowledgements}
The work of J. H. Jung and J. O. Lee was partially supported by National Research Foundation of Korea under grant number NRF-2019R1A5A1028324. The work of H. W. Chung was partially supported by National Research Foundation of Korea under grant number 2017R1E1A1A01076340 and by the Ministry of Science and ICT, Korea, under an ITRC Program, IITP-2019-2018-0-01402.

\appendix

\section{Simulations} \label{sec:simulation}
In this section, we numerically check our main results, the improvement of PCA and the error of the proposed test for hypothesis testing.
\subsection{Improved PCA with Entrywise Transformation}
We consider the data with non-Gaussian noise. We let the density function of the noise be the bimodal distribution with unit variance, defined as
\[
g(x)=\frac{1}{\sqrt{2\pi}}\left(e^{-2(x-\sqrt{3}/2)^2}+e^{-2(x+\sqrt{3}/2)^2}\right),
\]
which is a convolution of Gaussian and Rademacher random variables; more precisely, $g$ is the density function of a random variable
\[
\frac{1}{2} \caN + \frac{\sqrt{3}}{2}\caR,
\]
where $\caN$ is a standard Gaussian random variable and $\caR$ is a Rademacher random variable, independent to each other.

We sample $W_{ij}$ independently from the density $g$ and let $X_{ij}=W_{ij}/\sqrt{N}$.
We let $\bsu = (u_1, u_2, \dots, u_M)^T$ and $\bsv = (v_1, v_2, \dots, v_N)^T$, where $\sqrt{M} u_i$'s and $\sqrt{M} v_j$'s are i.i.d. Rademacher random variables for $i=1, 2, \dots, M$ and $j=1, 2, \dots, N$. The data matrix with additive spike is 
\beq \label{eq:generated_data}
Y= \sqrt{\lambda} \bsu \bsv^T + X.
\eeq
When we apply the entrywise transformation, defined in~\eqref{eq:transformed-cov}, with $\alpha=0$ to the data matrix, we get
\beq
\wt Y_{ij} = \frac{1}{\sqrt{\fg N}} h(\sqrt{N} Y_{ij})
\eeq
where 
\beq
h(x)=-\frac{g'(x)}{g(x)}=\frac{2\left(\sqrt{3}-e^{4\sqrt{3}x\left(\sqrt{3}-2x\right)}+2x\right)}{1+e^{4\sqrt{3}x}}
\eeq
and $\fg = \int_{-\infty}^{\infty} \frac{(g'(x))^2}{g(x)} \dd x\approx 2.50810$.
The size of the data matrix is set to be $M=1024$, $N=2048$, and the ratio $\rat=M/N=1/2$.
From the Marchenko-Pastur law, the transition of the largest eigenvalue occurs at $\sqrt{\rat}\approx 0.7071$.
After the transformation, on the other hand, from Theorem~\ref{thm:trans-mean} it is expected that the transition of the largest eigenvalue occurs at
$
\frac{\sqrt{\rat}}{\fg} \approx 0.2819.
$
We set the SNR 
\beq\label{eq:lam_sim}
\lambda=\frac{\sqrt{\rat}+\frac{\sqrt{\rat}}{\fg}}{2}
\eeq
to observe the transitions of the largest eigenvalue after the transformation.  

\begin{figure*}
	\centering
	\begin{minipage}[b]{.45\textwidth}
		\vskip 0.2in
		\begin{center}
			\centerline{\includegraphics[width=\columnwidth]{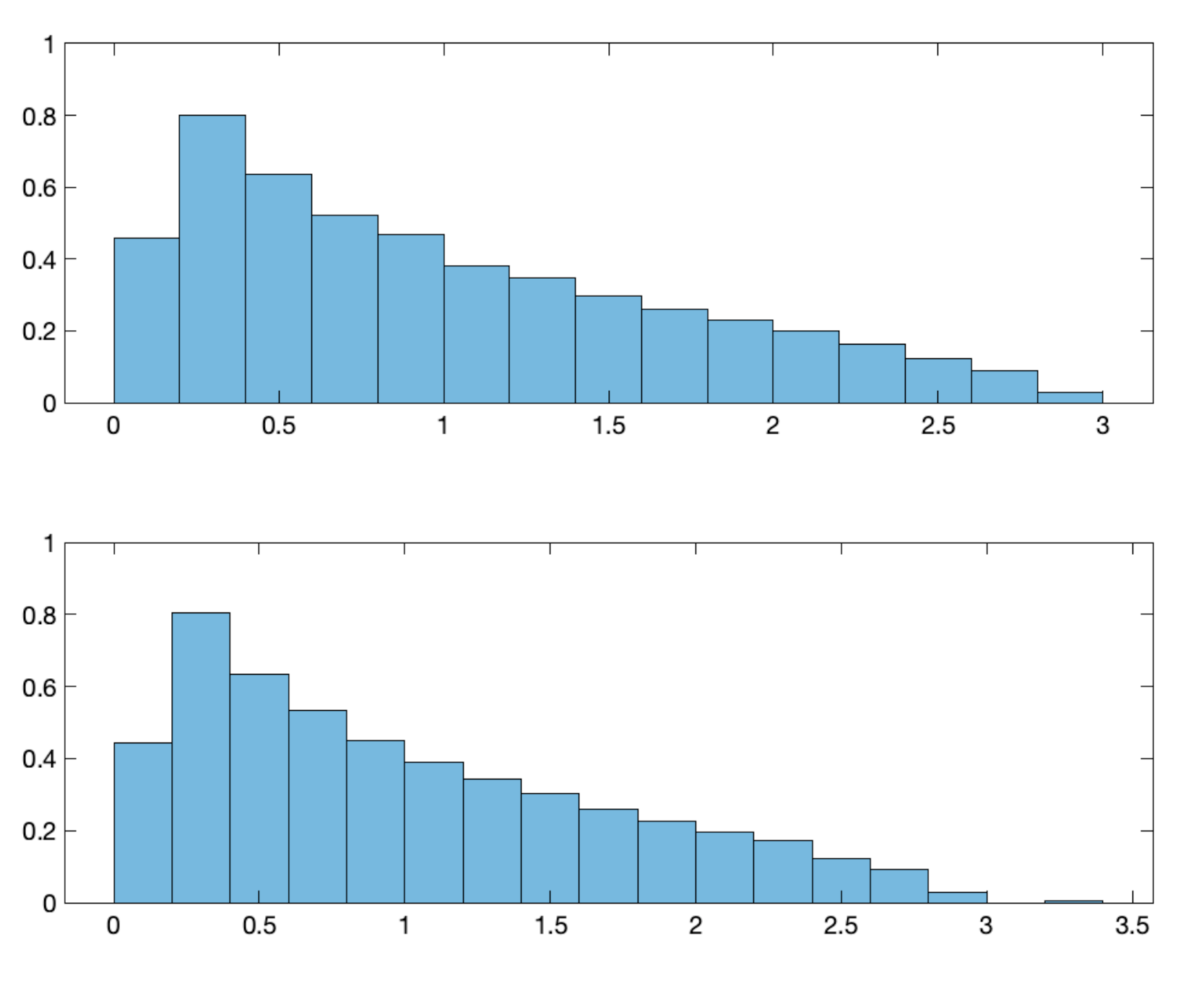}}
			\caption{The spectrum of the sample covariance matrix ($M=1024, N=2048, \lambda=0.4945$) with bimodal noise, before (above) and after (below) the entrywise transformation. An isolated eigenvalue can be seen only after the entrywise transformation. }
			\label{fig:spectrum}
		\end{center}
		\vskip -0.2in
		
	\end{minipage}\qquad
	\begin{minipage}[b]{.45\textwidth}
		\vskip 0.2in
		\begin{center}
			\centerline{\includegraphics[width=\columnwidth]{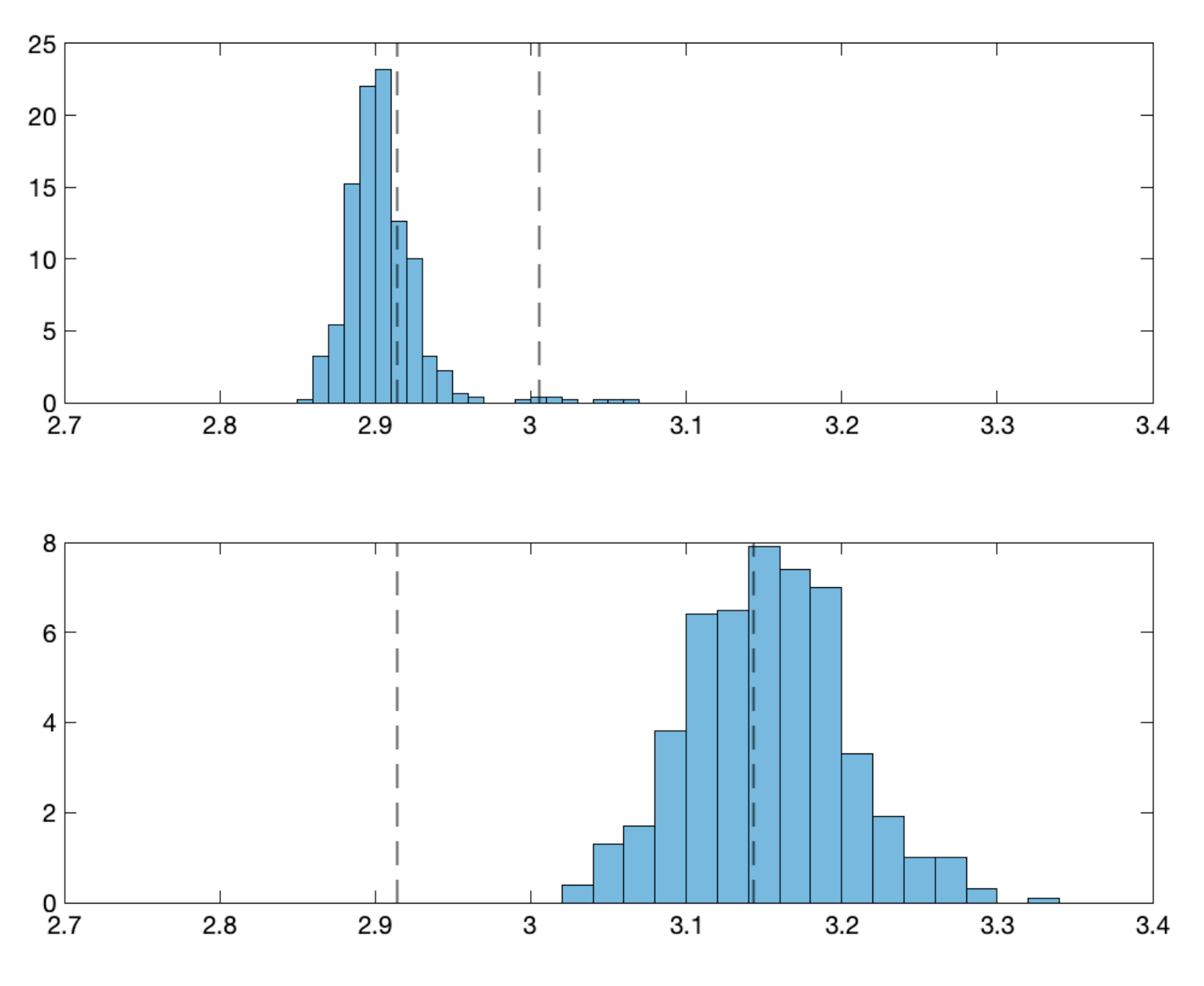}}
			\caption{Histograms of the largest eigenvalues of $YY^T$ (top) and $\wt Y \wt Y^T$ (bottom) from Monte Carlo simulations over 500 trials, with vertical lines indicating the theoretical results. The entrywise transformation makes the largest eigenvalue separated from the spectrum of the random matrix. }
			\label{fig:largest}
		\end{center}
		\vskip -0.2in
	\end{minipage}
\end{figure*}

In Fig.~\ref{fig:spectrum}, we compare the spectrum of the sample covariance matrices, $YY^T$ (top) and $\wt Y \wt Y^T$ (bottom).
An isolated eigenvalue can be seen only in the spectrum in the bottom, which is the case after the entrywise transformation.

In Fig.~\ref{fig:largest}, we compare the histograms of the largest eigenvalues of $YY^T$ (top) and $\wt Y \wt Y^T$ (bottom) from Monte Carlo simulations over 500 trials, and compare them with the theoretical results (vertical lines), $(1+\sqrt{\rat})^2$ and $(1+\lambda)(1+\frac{\rat}{\lambda})$ for $YY^T$ (top) and $(1+\sqrt{\rat})^2$ and $(1+\lambda \fg)(1+\frac{\rat}{\lambda \fg})$ for $\wt Y \wt Y^T$ (bottom).
From the simulations, it can be checked that at the same value of $\lambda$ (in~\eqref{eq:lam_sim}) PCA works only for the bottom case (after the entrywise transformation), and the largest eigenvalue closely matches the theoretical result.

\subsection{Hypothesis Testing with LSS estimator} \label{subsec:simul_LSS}

We consider the case where the noise matrix $X$ is Gaussian and the signal $\bsu = (u_1, u_2, \dots, u_M)^T$ and $\bsv = (v_1, v_2, \dots, v_N)^T$, where $\sqrt{M} u_i$'s and $\sqrt{N} v_j$'s are i.i.d. Rademacher random variables for $i=1, 2, \dots, M$ and $j=1, 2, \dots, N$. Let the data matrix $Y = \sqrt{\lambda} \bsu \bsv^T + X$. The parameters are $w_2 = 2$ and $w_4 = 3$. 

In the numerical simulation done in Matlab, we generated 10,000 independent samples of the $256 \times 512$ data matrix $Y$ under $\bsH_0$ (without signal $\lambda=0$) and $\bsH_1$ (with signal $\lambda=\SNR$), respectively, varying SNR $\SNR$ from $0$ to $0.5$. To apply Algorithm \ref{alg:test} proposed in Section \ref{subsec:ht}, we computed
\beq \begin{split}
	L_{\SNR} = -\log \det \left( \left(1+\frac{\rat}{\SNR} \right)(1+\SNR)I - YY^T \right) 
	+M \left[\frac{\SNR}{\rat} - \log\left(\frac{\SNR}{\rat}\right) -\frac{1-\rat}{\rat}\log(1+\SNR) \right],
\end{split} \eeq
accept $\bsH_0$ if $L_{\SNR} \leq -\log\left(1-\frac{\SNR^2}{\rat}\right)$ and reject $\bsH_0$ otherwise. The limiting error of the test is
\beq \label{eq:limit_error_1a}
\erfc \left( \frac{\sqrt{V_0}}{4\sqrt 2} \right) = \erfc \left( \frac{1}{4} \sqrt{-\log \left( 1- \frac{\SNR^2}{\rat} \right)} \right),
\eeq
where $V_0$ is the variance in \eqref{eq:var_test} and $\erfc(\cdot)$ is the complementary error function.

\begin{figure}[t]
	\vskip 0.2in
	\begin{center}
		\centerline{\includegraphics[width=0.75\columnwidth]{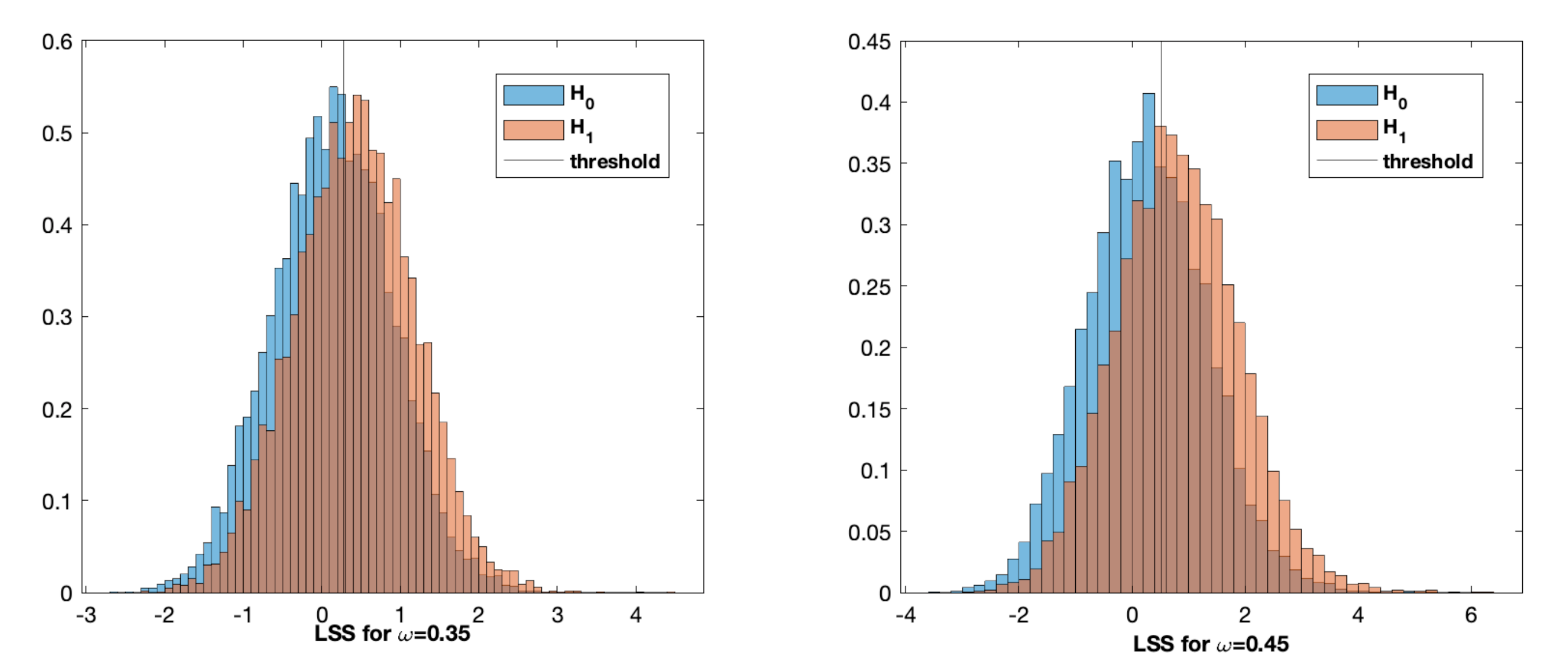}}
		\caption{The histograms of the test statistic $L_\omega$ under $H_0$ and $H_1$, respectively, for the Gaussian noise with $\omega=0.35$ and $\omega=0.45$. The threshold for the test is given by $-\log(1-\omega^2/d)$.}
		\label{fig:histo_HT}
	\end{center}
	\vskip -0.2in
\end{figure}

In Figure \ref{fig:histo_HT}, we plot the histograms of the test statistic $L_{\SNR}$ under $\bsH_0$ and under $\bsH_1$, respectively, with the test threshold $-\log\left(1-\frac{\SNR^2}{\rat}\right)$ for $\SNR = 0.35$ and $\SNR = 0.45$. It can be shown that the difference of the means of $L_{\SNR}$ under $\bsH_0$ and under $\bsH_1$ is larger for $\SNR = 0.45$.

\section{Proof of Theorems for improved PCA} \label{sec:proof_trans}

In this section, we prove our main results on the improved PCA, Theorems \ref{thm:trans-mean} and \ref{thm:trans-cov}. 

\subsection{Preliminaries}
We first introduce the following notions, which provide a simple way of making precise statements regarding the bound up to small powers of $N$ that holds with probability higher than $1-N^{-D}$ for all $D>0$.

\begin{defn}[Overwhelming probability]\label{def:overwhelming}
	We say that an event (or family of events) $\Omega$ holds with overwhelming probability if for all (large) $D>0$ we have $\P(\Omega) \leq N^{-D}$ for any sufficiently large $N$.
\end{defn}

\begin{defn}[Stochastic domination]\label{def:stocdom}
	Let
	\begin{equation*}
	\xi = \pb{\xi^{(N)}(u) \;:\; N \in \N, u \in U^{(N)}} \,, \qquad
	\zeta = \pb{\zeta^{(N)}(u) \;:\; N \in \N, u \in U^{(N)}}
	\end{equation*}
	be two families of random variables, where $U^{(N)}$ is a possibly $N$-dependent parameter set. 
	We say that $\xi$ is \emph{stochastically dominated by $\zeta$, uniformly in $u$,} if for all (small) $\epsilon > 0$ and (large) $D > 0$
	\begin{equation*}
	\sup_{u \in U^{(N)}} \P \left({|\xi^{(N)}(u)| > N^\epsilon \zeta^{(N)}(u)}\right) \;\leq\; N^{-D}
	\end{equation*}
	for any sufficiently large $N\geq N_0(\varepsilon, D)$.
	Throughout this appendix, the stochastic domination will always be uniform in all parameters, including matrix indices and the spectral parameter $z$.
	
	We write $\xi \prec \zeta$ or $\xi = \caO_\prec(\zeta)$, if $\xi$ is stochastically dominated by $\zeta$, uniformly in $u$.
\end{defn}

Under Assumption \ref{assump:entry}, we have
\[
\| \bsu \| = 1+\caO_\prec(N^{-\phi}),
\]
which can be proved by applying the Markov inequality with the bounded moment assumption on $\bsu$; more precisely, since $u_k$'s are independent, for any (large) $p$,
\beq \begin{split}
	\p \left(\left|\left(\sum_{k=1}^M u_k^2\right)-1\right|\geq N^\varepsilon N^{-\phi}\right) \leq\frac{\E\left|\sum_{k=1}^M\left(u_k^2-\frac1{M}\right)\right|^{2p}}{N^{2p\varepsilon} N^{-2p\phi}} \leq C\frac{N^{-2p\phi}}{N^{p\varepsilon} N^{-2p\phi}} = CN^{-p\varepsilon}
\end{split} \eeq
for some constant $C$. Similarly, $\| \bsv \| = 1+\caO_\prec(N^{-\phi})$.

We will use the following result for the resolvents, which is called an isotropic Marchenko--Pastur law.
\begin{lem}[Isotropic local Marchenko--Pastur law] \label{lem:local_MP}
	Suppose that $z \in \R$ outside an open interval containing $[d_-, d_+]$. Let $s(z)$ be the Stieltjes transform of the Marchenko--Pastur law, which is also given by
	\beq
	s(z) = \frac{(1-\rat-z)+ \sqrt{(1-\rat-z)^2 -4\rat z}}{2\rat z}.
	\eeq
	Then,
	\[
	\langle \bsv, (X^T X-zI)^{-1} \bsv \rangle = -\left( \frac{1}{zs(z)}+1 \right) \| \bsv \|^2 + \caO_\prec(N^{-\frac{1}{2}})
	\]
	and
	\[	
	\langle X^T \bsu, (X^T X-zI)^{-1} X^T\bsu\rangle = (zs(z)+1) \| \bsu \|^2 +\caO_\prec(N^{-\frac{1}{2}}).
	\]
\end{lem}
See Theorem 2.5 of \cite{Bloemendal-Erdos-Knowles-Yau-Yin2014} (also Lemma 3.7 of \cite{BKYY2016}) for the proof of Lemma \ref{lem:local_MP}.

The following concentration inequality will be frequently used in the proof, which is sometimes called the large deviation estimate in random matrix theory.
\begin{lem}[Large deviation estimate] \label{lemma: LDE}
	Let $\pb{\xi_i^{(N)}}$ and $\pb{\zeta_i^{(N)}}$ be independent families of random variables and 
	$\pb{a_{ij}^{(N)}}$ and $\pb{b_i^{(N)}}$ be deterministic; here $N \in \N$ and $i,j = 1, \dots, N$. Suppose that complex-valued random variables $\xi_i^{(N)}$ and $\zeta_i^{(N)}$ are independent and satisfy for all $p\geq 2$ that
	\beq\label{cond on X}
	\E \xi \;=\; 0\,, \qquad \E \abs{\xi}^p \;\leq\; \frac{C_p}{NB^{p-2}}
	\eeq
	for some $B \leq N^{1/2}$ and some ($N$-independent) constant $C_p$. Then we have the bounds
	\begin{align} \label{LDE}
	\sum_i b_i \xi_i &\;\prec\; \pbb{\frac1{N}\sum_i \abs{b_i}^2}^{1/2}+\frac{\max_i|b_i|}{B}\,,
	\\ \label{two-set LDE}
	\sum_{i,j} a_{ij} \xi_i \zeta_j &\;\prec\; \pbb{\frac1{N^2}\sum_{i\neq j} \abs{a_{ij}}^2}^{1/2}+\frac{\max_{i\neq j}|a_{ij}|}{B}+\frac{\max_{i}|a_{ii}|}{B^2}\,,
	\\ \label{offdiag LDE}
	\sum_{i \neq j} a_{ij} \xi_i \xi_j &\;\prec\; \pbb{\frac1{N^2}\sum_{i \neq j} \abs{a_{ij}}^2}^{1/2}+\frac{\max_{i\neq j}|a_{ij}|}{B}\,.
	\end{align}
	If the coefficients $a_{ij}^{(N)}$ and $b_i^{(N)}$ depend on an additional parameter $u$, then all of these estimates are uniform in $u$, i.e. $N_0= N_0(\varepsilon, D)$ in the definition of $\prec$ depends not on $u$ but only on the constant $C$ from \eqref{cond on X}. 
	
	If $B = N^{1/2}$, the bounds can further be simplified to
	\beq
	\sum_i b_i \xi_i \prec \pbb{\frac1{N}\sum_i \abs{b_i}^2}^{1/2}, \quad \sum_{i,j} a_{ij} \xi_i \zeta_j \prec \pbb{\frac1{N^2}\sum_{i, j} \abs{a_{ij}}^2}^{1/2}, \quad
	\sum_{i \neq j} a_{ij} \xi_i \xi_j \prec\pbb{\frac1{N^2}\sum_{i \neq j} \abs{a_{ij}}^2}^{1/2}.
	\eeq
\end{lem}

\begin{proof}
	These estimates are an immediate consequence of Lemma 3.8 in \cite{Erdos-Knowles-Yau-Yin2013a}.
\end{proof}

\subsection{Proof of Theorem \ref{thm:trans-mean}} \label{sec:proof_trans-mean}

We first prove the behavior of the largest eigenvalue described in Section \ref{subsec:BBP}, which we will call the BBP result, in our setting, following the strategy of \cite{Raj2011,benaych2012singular}. Note that the largest eigenvalue of $YY^T$ is equal to the largest eigenvalue of $Y^T Y$. Consider the identity
\beq \begin{split}
	Y^T Y -zI &= (X+ \lambda \bsu \bsv^T)^T (X+ \lambda \bsu \bsv^T) - zI \\
	&= (X^T X-zI)(I + (X^T X-zI)^{-1} (\lambda X^T \bsu \bsv^T + \lambda \bsv \bsu^T X + \lambda^2 \| \bsu \|^2 \bsv \bsv^T)).
\end{split} \eeq
Thus, if $z$ is an eigenvalue of $YY^T$ but not of $XX^T$, then it satisfies
\[
\det (I + (X^T X-zI)^{-1} (\lambda X^T \bsu \bsv^T + \lambda \bsv \bsu^T X + \lambda^2 \| \bsu \|^2 \bsv \bsv^T)) = 0,
\]
which also implies that $-1$ is an eigenvalue of 
\[
T \equiv T(z) := (X^T X-zI)^{-1} (\lambda X^T \bsu \bsv^T + \lambda \bsv \bsu^T X + \lambda^2 \| \bsu \|^2 \bsv \bsv^T).
\]
The rank of $T$ is at most $2$, with
\[ \begin{split}
T \cdot (X^T X-zI)^{-1} X^T \bsu &= \lambda \langle \bsv, (X^T X-zI)^{-1} X^T \bsu \rangle (X^T X-zI)^{-1} X^T \bsu \\
& \qquad + \lambda \langle X^T \bsu, (X^T X-zI)^{-1} X^T \bsu \rangle (X^T X-zI)^{-1} \bsv \\
& \qquad + \lambda^2 \| \bsu \|^2 \langle \bsv, (X^T X-zI)^{-1} X^T \bsu \rangle (X^T X-zI)^{-1} \bsv
\end{split} \]
and
\[ \begin{split}
T \cdot (X^T X-zI)^{-1} \bsv &= \lambda \langle \bsv, (X^T X-zI)^{-1} \bsv \rangle (X^T X-zI)^{-1} X^T \bsu \\
& \qquad + \lambda \langle X^T \bsu, (X^T X-zI)^{-1} \bsv \rangle (X^T X-zI)^{-1} \bsv \\
& \qquad + \lambda^2 \| \bsu \|^2 \langle \bsv, (X^T X-zI)^{-1} \bsv \rangle (X^T X-zI)^{-1} \bsv.
\end{split} \]
In particular, an eigenvector of $T$ is a linear combination of $(X^T X-zI)^{-1} X^T \bsu$ and $(X^T X-zI)^{-1} \bsv$.

Suppose that $a (X^T X-zI)^{-1} X^T \bsu + b (X^T X-zI)^{-1} \bsv$ is an eigenvector of $T$ with the corresponding eigenvalue $-1$.  Since $\| X \|, \| (X^T X-zI)^{-1} \| \prec 1$, from Lemma \ref{lemma: LDE}, 
\[ \begin{split}
&\langle \bsv, (X^T X-zI)^{-1} X^T \bsu \rangle = \sum_{i, j} \left( (X^T X-zI)^{-1} X^T \right)_{ij} v_i u_j \\
&\prec \left( \frac{1}{N^2} \sum_{i \neq j} \left| \left( (X^T X-zI)^{-1} X^T \right)_{ij} \right|^2 \right)^{1/2} + N^{-\phi} \max_{i, j} \left| \left( (X^T X-zI)^{-1} X^T \right)_{ij} \right| \\
&\prec \left( \frac{1}{N} \| (X^T X-zI)^{-1} X^T \|^2 \right)^{1/2} + N^{-\phi} \| (X^T X-zI)^{-1} X^T \| \prec N^{-\phi}.
\end{split} \]
Thus, from Lemma \ref{lem:local_MP},
\beq \begin{split} \label{eq:ab_equation}
	&-\left( a (X^T X-zI)^{-1} X^T \bsu + b (X^T X-zI)^{-1} \bsv \right) = T \left( a (X^T X-zI)^{-1} X^T \bsu + b (X^T X-zI)^{-1} \bsv \right) \\
	&= -b\lambda \left( \frac{1}{zs(z)}+1 \right) (X^T X-zI)^{-1} X^T \bsu \\
	&\qquad + a\lambda (zs(z)+1) (X^T X-zI)^{-1} \bsv - b\lambda^2 \left( \frac{1}{zs(z)}+1 \right) (X^T X-zI)^{-1} \bsv + \bsth
\end{split} \eeq
for some $\bsth$, which is a linear combination of $(X^T X-zI)^{-1} X^T \bsu$ and $(X^T X-zI)^{-1} \bsv$ with $\| \bsth \| = \caO_{\prec}(N^{-\phi})$.

Since $\bsu$, $\bsv$, and $X$ are independent, $(X^T X-zI)^{-1} X^T \bsu$ and $(X^T X-zI)^{-1} \bsv$ are linearly independent with overwhelming probability. Thus, from \eqref{eq:ab_equation}, 
\[ \begin{split}
-a &= -b\lambda \left( \frac{1}{zs(z)}+1 \right) + \caO_{\prec}(N^{-\phi}), \\
-b &= a\lambda (zs(z)+1) - b\lambda^2 \left( \frac{1}{zs(z)}+1 \right) + \caO_{\prec}(N^{-\phi}).
\end{split} \]
It is then elementary to check that
\[
\lambda^2 (zs(z)+1) +1 = \caO_{\prec}(N^{-\phi}),
\]
which has the solution
\[
z = (1+\lambda)(1+ \frac{\rat}{\lambda}) + \caO_{\prec}(N^{-\phi})
\]
if and only if $\lambda > \sqrt{\rat}$. This proves the BBP result in our setting.

We now turn to the proof of Theorem \ref{thm:trans-mean}. For the i.i.d. prior, suppose that a function $\f$ and its all derivatives are polynomially bounded in the sense of Assumption \ref{assump:entry}. Following the proof of Theorem 4.8 in \cite{Perry2018}, we define the error term from the local linear estimation of $\f(\sqrt{N} Y_{ij})$ by
\[
\f(\sqrt{N}Y_{ij})=\f(\sqrt{N}X_{ij})+\sqrt{\lambda N} u_i v_j \f'(\sqrt{N}X_{ij})+\mathcal{R}_{ij}
\]
where 
\[
\mathcal{R}_{ij}=\frac{1}{2} \f''(\sqrt{N}X_{ij}+e_{ij}) \lambda N u^2_i v^2_j
\] 
for some $|e_{ij}|\leq|\sqrt{\lambda N}u_i v_j|$. The Frobenius norm of $\mathcal{R}$ is bounded as
\[
\|\mathcal{R}\|^2_F= \Tr \mathcal{R}^T \mathcal{R}=\frac{\lambda^2 N^2}{4}\sum_{i=1}^{M}\sum_{j=1}^{N}u_i^4v_j^4 \f''(\sqrt{N}X_{ij}+e_{ij})^2\leq\frac{\lambda^2 N^{2-4\phi}}{4}\sum_{i=1}^{M}\sum_{j=1}^{N}u_i^2v_j^2 \f''(\sqrt{N}X_{ij}+e_{ij})^2.
\]
Since $\f''$ is polynomially bounded, $\f''(\sqrt{N}X_{ij}+e_{ij})$ is uniformly bounded by an $N$-independent constant. Thus, with overwhelming probability,
\[
\|\mathcal{R}\|^2\leq\|\mathcal{R}\|^2_F\leq C\lambda^2N^{2-4\phi}.
\]

Next, we approximate $\f(\sqrt{N}X_{ij})$ by its mean. Let 
\[ \label{eq:local_linear}
\mathcal{E}_{ij}=\f'(\sqrt{N}X_{ij})-\E [\f'(\sqrt{N}X_{ij})], \qquad \Delta_{ij}=\sqrt{\lambda N}u_i v_j \mathcal{E}_{ij}.
\]
Then, $\|\Delta\|\prec N^{\frac12-2\phi}\|\mathcal{E}\|$ and, since the entries of matrix $\mathcal{E}$ are i.i.d., centered and with finite moments, its norm $\|\mathcal{E}\|=O(\sqrt{N})$ with overwhelming probability. (See, e.g., \cite{BKYY2016}.) Thus, $\| \Delta \| = \caO_{\prec}(N^{1-2\phi})$.

Set
\[
\mf := \E [\f'(\sqrt{N}X_{ij})], \quad \vf := \E [\f(\sqrt{N}X_{ij})^2], \quad \widehat\lambda := \lambda \mf^2/\vf,
\]
and
\[
\F_{ij} := \frac{1}{\sqrt{N\vf}}\f(\sqrt{N}X_{ij}).
\]
We have proved so far that the difference of the largest eigenvalue of $\F + \widehat\lambda^{\frac{1}{2}} \bsu \bsv^T$ and that of the matrix
\[
\left( \frac{1}{\sqrt{N\vf}}\f(\sqrt{N}Y_{ij}) \right)
\]
is $\caO_{\prec}(N^{\frac12-2\phi})$, which is $o(1)$ with overwhelming probability for $\phi > \frac{1}{4}$. The BBP result holds the matrix $\F + \widehat\lambda^{\frac{1}{2}} \bsu \bsv^T$, which is another (additive) spiked rectangular matrix. This shows that the BBP result also holds for $\wt Y$ with SNR $\widehat\lambda := \lambda \mf^2/\vf$. This proves Theorem \ref{thm:trans-mean} for the i.i.d. prior.

For the spherical prior, we replace $\bsu$ and $\bsv$ by normalized Gaussian vectors. More precisely, let $\bsw \in \R^M$ and $\wt\bsw \in \R^N$ be random vectors whose entries are i.i.d. standard Gaussian. From the spherical symmetry of Gaussian, we find that $\bsu$ and $\bsv$ have the same distributions as $\bsw/\| \bsw \|$ and $\wt\bsw / \| \wt \bsw \|$. If we replace $\bsu$ and $\bsv$ by $\bsw/\sqrt{M}$ and $\wt\bsw/\sqrt{N}$, the spike prior is i.i.d. and the desired theorem holds. Since $\| \bsw \| = M + o(M)$ and $\| \wt \bsw \| = N + o(N)$ with overwhelming probability, the change of $\lambda$ due to the replacement is negligible in the limit $M, N \to \infty$. This proves Theorem \ref{thm:trans-mean} for the spherical prior.

\subsection{Proof of Theorem \ref{thm:trans-cov}}\label{sec:proof_trans-cov}

We assume the i.i.d. prior; the proof for the spherical prior easily follows from the result with the i.i.d. prior as in the proof of Theorem \ref{thm:trans-mean}. As in the additive case, we assume that a function $\f$ and its all derivatives are polynomially bounded and consider the local linear approximation of $\f(\sqrt{N} Y_{ij})$,
\beq \label{eq:local_linear2}
\f(\sqrt{N} Y_{ij} ) = \f(\sqrt{N} X_{ij}) + \gamma \sqrt{N} \E[\f'(\sqrt{N} X_{ij})] \sum_k u_i u_k X_{kj} + \mathcal{R}_{ij} + \gamma \Delta_{ij},
\eeq
where
\[
\mathcal{R}_{ij}=\frac{1}{2}\f''\Big(\sqrt{N} X_{ij} + \theta\gamma\sum_k u_i u_k \sqrt{N} X_{kj} \Big)\left(\gamma  \sum_k u_i u_k \sqrt{N}X_{kj} \right)^2
\]
for some $\theta \in [-1, 1]$ and 
\[
\Delta_{ij}=\sqrt{N} \caE_{ij} \sum_k u_i u_k X_{kj}, \quad \caE_{ij}=\f'(\sqrt{N} X_{ij})-\E[\f'(\sqrt{N} X_{ij})].
\]

For any unit vectors $\bsa = (a_1, a_2, \dots, a_M)$ and $\bsb = (b_1, b_2, \dots, b_N)$,
\[
\bsa^T\Delta \bsb=\sum_{i,j}a_iu_i\caE_{ij}b_j\left(\sum_k u_k \sqrt{N}X_{kj}\right) 
=\sum_{i,j}a_iu_i^2b_j\caE_{ij}\sqrt{N}X_{ij}+\sum_{i,j}a_iu_i\caE_{ij}b_j\left(\sum_{k\neq i} u_k \sqrt{N}X_{kj}\right)
\]
From the concentration inequalities such as Lemma \ref{lemma: LDE},
\beq \label{eq:u_concentration}
\sum_k u_i u_k \sqrt{ N}X_{kj}\prec |u_i|\left(\sum_ku_k^2\right)^{1/2}\prec N^{-\phi}.
\eeq
Recall that $\| \caE \| = O(\sqrt{N})$ with overwhelming probability. Note that, by Assumption \ref{assump:entry}, the density function $\f$ have to be
an odd function. Further, since $\f$ is an odd function (hence $x\f'(x)$ is an odd function of $x$), the norm of the matrix whose $(i, j)$-entry is $\caE_{ij} \sqrt{N} X_{ij}$ is also $O(\sqrt{N})$. 
Thus,
\[
\bsa^T\Delta \bsb \prec N^{-2\phi} + N^{\frac{1}{2}-\phi},
\]
which shows that $\| \Delta \| \prec N^{\frac{1}{2}-\phi}$. Moreover, since $\f''$ is polynomially bounded, following the proof of Theorem \ref{thm:trans-mean} with \eqref{eq:u_concentration},
\[
\|\mathcal{R} \|^2 \leq \| \mathcal{R}\| _F^2 \leq CN^{2-4\phi}.
\]
Thus, as in the additive case, the error terms $\mathcal{R}_{ij}$ and $\Delta_{ij}$ in \eqref{eq:local_linear2} are negligible when finding the limit of the largest eigenvalue of the transformed matrix.

Set
\[
\mf := \E [\f'(\sqrt{N}X_{ij})], \quad \vf := \E [\f(\sqrt{N}X_{ij})^2], \quad \ef = \E [\sqrt{N}X_{ij} \f(\sqrt{N}X_{ij})], \quad \widehat\gamma :=\gamma \mf/\sqrt{\vf},
\]
and
\[
\F_{ij} := \frac{1}{\sqrt{N\vf}}\f(\sqrt{N}X_{ij}).
\]
With the approximation \eqref{eq:local_linear2}, we now focus on the largest eigenvalue of
\[
(\F+\wh\gamma\bsu\bsu^TX)^T(\F+\wh\gamma\bsu\bsu^TX).
\]
Note that the assumption on the polynomial boundedness of $\f$ implies that the matrix $\F$ is also a rectangular matrix satisfying the assumptions in Definition \ref{defn:rect}.

Let $G(z)$ and $\caG(z)$ be the resolvents
\[
G \equiv G(z) := (\F\F^T -zI)^{-1}, \qquad \caG \equiv \caG(z) := (\F^T \F -zI)^{-1}
\]
for $z \in \R$ outside an open interval containing $[d_-, d_+]$. We note that the following identities hold for $G(z)$ and $\caG(z)$:
\beq\label{eq:relation between g}
G(z)\F = \F \caG(z), \qquad \F^T G(z) \F=I+z\caG(z).
\eeq

As in the proof of Theorem \ref{thm:trans-mean}, we consider
\beq\begin{split} \label{eq:determinant}
	&(\F+\wh\gamma\bsu\bsu^TX)^T(\F+\wh\gamma\bsu\bsu^TX)-zI 
	\\&=(\F^T\F-zI)(I+(\F^T\F-zI)^{-1}(\wh\gamma X^T\bsu\bsu^T\F+\wh\gamma \F^T \bsu\bsu^TX+ \|\bsu \|^2 \wh\gamma^2X^T\bsu\bsu^TX)).
\end{split}\eeq
Let 
\[
L \equiv L(z)=\caG(z)(\wh\gamma X^T\bsu\bsu^T\F+\wh\gamma \F^T \bsu\bsu^TX+ \|\bsu \|^2 \wh\gamma^2X^T\bsu\bsu^TX),
\]
Then, as in the proof of Theorem \ref{thm:trans-mean}, if $z$ is an eigenvalue of $(\F+\wh\gamma\bsu\bsu^T X)^T (\F+\wh\gamma\bsu\bsu^T X)$ (but not of $\F^T \F$), $-1$ is an eigenvalue of $L(z)$. Again, the rank of $L$ is at most $2$, with
\beq \begin{split} \label{eq:eigenvector}
	L \cdot \caG \F^T\bsu &=\wh\gamma \langle\bsu,\F\caG \F^T\bsu\rangle \cdot \caG X^T\bsu+\wh\gamma \langle\bsu,X\caG \F^T\bsu\rangle \cdot \caG \F^T\bsu+\|\bsu\|^2\wh\gamma^2\langle\bsu,X\caG \F^T\bsu\rangle \cdot \caG X^T\bsu, \\
	L \cdot \caG X^T\bsu &=\wh\gamma \langle\bsu,\F\caG X^T\bsu\rangle \cdot \caG X^T\bsu+\wh\gamma \langle\bsu, X\caG X^T\bsu\rangle \cdot \caG \F^T\bsu+\|\bsu\|^2\wh\gamma^2\langle\bsu,X\caG X^T\bsu\rangle \cdot \caG X^T\bsu,
\end{split} \eeq
and an eigenvector of $L$ is a linear combination of $\caG \F^T\bsu$ and $\caG X^T\bsu$.

In the simplest case where $\f$ is the identity mapping, $\F=X$, hence the rank of $L$ is $1$, and the eigenvalue equation \eqref{eq:eigenvector} is simplified to
\beq \label{eq:eigenvector_i}
L \cdot \caG \F^T\bsu =\wh\gamma\langle\bsu,\F\caG \F^T\bsu\rangle \cdot \caG \F^T\bsu+\wh\gamma\langle\bsu, \F\caG \F^T\bsu\rangle \cdot \caG \F^T\bsu+\|\bsu\|^2\wh\gamma^2\langle\bsu,\F\caG \F^T\bsu\rangle \cdot \caG \F^T\bsu.
\eeq
In this case, $\caG \F^T\bsu$ is an eigenvector of $L$ corresponding to the eigenvalue $-1$, i.e., $L \cdot \caG \F^T\bsu = -\caG \F^T\bsu$. The right side of \eqref{eq:eigenvector_i} can be approximated as follows, which is a direct consequence of the isotropic local Marchenko--Pastur law (e.g., Theorem 2.5 of \cite{Bloemendal-Erdos-Knowles-Yau-Yin2014}).

With the isotropic local Marchenko--Pastur law, \eqref{eq:eigenvector_i} can be approximated by a deterministic equation on $z$ (and $s(z)$), and the location of the largest eigenvalue can be proved by solving the equation. In a general case where $\F$ is not a multiple of $X$ and the vectors $\caG \F^T\bsu$ and $\caG X^T\bsu$ are linearly independent, however, the eigenvalue equation \eqref{eq:eigenvector} contains other terms $\langle\bsu,X\caG \F^T\bsu\rangle$, $\langle\bsu,\F\caG X^T\bsu\rangle$, and $\langle\bsu,X\caG X^T\bsu\rangle$, which cannot be estimated by Lemma \ref{lem:local_MP}. For these terms, we use the following lemma.
\begin{lem} \label{lem:local_MP2}
	Suppose that the assumptions in Lemma \ref{lem:local_MP} hold. Then,
	\[
	\langle\bsu,X\caG \F^T\bsu\rangle = \langle\bsu,\F\caG X^T\bsu\rangle = \frac{\ef}{\sqrt{\vf}}(zs(z)+1) +\caO_\prec(N^{-\phi})
	\]
	and
	\[
	\langle\bsu,X\caG X^T\bsu\rangle = \frac{\ef^2}{\vf} zs(z) \left( \rat s(z) + \frac{\rat-1}{z} \right)^2 + \rat s(z) + \frac{\rat-1}{z} +\caO_\prec(N^{-\phi}).
	\]
\end{lem}
We defer the proof to Appendix \ref{subsec:proof_lem}.

With Lemma \ref{lem:local_MP2}, we are ready to finish the proof. From the definition of $s(z)$ in Lemma \ref{lem:local_MP}, we notice that
\beq
s(z)=\frac{1}{1-\rat - \rat z s(z)-z},
\eeq
or
\beq
z\left( \rat s(z) + \frac{\rat-1}{z} \right) = -\frac{1}{s(z)} -z.
\eeq
Set $\sigma(z) := zs(z) + 1$.
From \eqref{eq:eigenvector_i},
\beq \begin{split} \label{eq:eigen_F}
	L \cdot \caG \F^T\bsu &=\wh\gamma \langle\bsu,\F\caG \F^T\bsu\rangle \cdot \caG X^T\bsu+\wh\gamma \langle\bsu,X\caG \F^T\bsu\rangle \cdot \caG \F^T\bsu+\|\bsu\|^2\wh\gamma^2\langle\bsu,X\caG \F^T\bsu\rangle \cdot \caG X^T\bsu \\
	&= \wh\gamma \sigma(z) \caG X^T\bsu+ \wh\gamma \sigma(z) \frac{\ef}{\sqrt{\vf}} \caG \F^T\bsu+ \wh\gamma^2 \frac{\ef}{\sqrt{\vf}} \sigma(z) \caG X^T\bsu+ \bsth_1 \,,
\end{split} \eeq
and
\beq \begin{split} \label{eq:eigen_X}
	L \cdot \caG X^T\bsu &=\wh\gamma \langle\bsu,\F\caG X^T\bsu\rangle \cdot \caG X^T\bsu+ \wh\gamma \langle\bsu, X\caG X^T\bsu\rangle \cdot \caG \F^T\bsu+\|\bsu\|^2\wh\gamma^2\langle\bsu,X\caG X^T\bsu\rangle \cdot \caG X^T\bsu \\
	&= \wh\gamma \sigma(z) \frac{\ef}{\sqrt{\vf}} \caG X^T\bsu+ \wh\gamma \left( \left(\sigma(z) + \frac{\sigma(z)}{\sigma(z)-1} \right) \frac{\ef^2}{\vf} - \frac{\sigma(z)}{\sigma(z)-1} \right) \caG \F^T\bsu \\
	&\qquad +\wh\gamma^2 \left( \left(\sigma(z) + \frac{\sigma(z)}{\sigma(z)-1} \right) \frac{\ef^2}{\vf} - \frac{\sigma(z)}{\sigma(z)-1} \right) \caG X^T\bsu + \bsth_2 \,,
\end{split} \eeq
for some $\bsth_1, \bsth_2$, which are linear combinations of $\caG \F^T\bsu$ and $\caG X^T\bsu$, with $\| \bsth_1 \|, \| \bsth_2 \| = \caO_\prec(N^{-\phi})$.

Suppose that $a\caG \F^T\bsu+b\caG X^T\bsu$ is an eigenvector of $L$ with the corresponding eigenvalue $-1$. From \eqref{eq:eigen_F}, \eqref{eq:eigen_X}, and the linear independence between $\caG \F^T\bsu$ and $\caG X^T\bsu$, we find the relation
\[ \begin{split}
-a &= a\wh\gamma \sigma(z) \frac{\ef}{\sqrt{\vf}} + \frac{b\wh\gamma \sigma(z)^2}{\sigma(z)-1} \cdot \frac{\ef^2}{\vf} - \frac{b\wh\gamma \sigma(z)}{\sigma(z)-1} + \caO(N^{-\phi}), \\
-b &= a\wh\gamma \sigma(z) + a \wh\gamma^2 \sigma(z) \frac{\ef}{\sqrt{\vf}} + b \wh\gamma \sigma(z) \frac{\ef}{\sqrt{\vf}} + \frac{b\wh\gamma^2 \sigma(z)^2}{\sigma(z)-1} \cdot \frac{\ef^2}{\vf} - \frac{b\wh\gamma^2 \sigma(z)}{\sigma(z)-1} + \caO(N^{-\phi}).
\end{split} \]
We then find that
\[
\frac{b}{a} \left( 1 +\wh\gamma \sigma(z) \frac{\ef}{\sqrt{\vf}} \right) + \wh\gamma \sigma(z) -\wh\gamma = \caO(N^{-\phi}),
\]
which implies that
\beq\label{res2}
1+2\wh\gamma \sigma(z) \frac{\ef}{\sqrt{\vf}}+ \wh\gamma^2 \sigma(z) =1+\left(\frac{2\gamma \mf\ef+\gamma^2\mf^2}{\vf}\right) \sigma(z) = \caO(N^{-\phi}).
\eeq 
From the explicit formula for $s$, it is not hard to check that \eqref{res2} holds if and only if
\[
\lambda_\f := \frac{2\gamma \mf\ef+\gamma^2\mf^2}{\vf} > \sqrt{\rat}
\]
and
\beq \label{eigen2}
z = (1+\lambda_\f) (1+\frac{\rat}{\lambda_\f}) +\caO(N^{-\phi}).
\eeq
Now, the desired theorem follows from the direct computation for the case $\f = h_{\alpha_g}$; see also Appendix \ref{subsec:variation_cov}.

\subsection{Optimal entrywise transformation} \label{sec:optimize_transform}

\subsubsection{Additive model} \label{subsec:variation_mean}
Recall that 
\[
\E [\f'(\sqrt{N}X_{ij})]=\mf, \qquad \E [\f(\sqrt{N}X_{ij})^2]=\vf.
\]
Following the proof of Theorem \ref{thm:trans-mean} in Appendix \ref{sec:proof_trans-mean}, it is not hard to see that the effective SNR is maximized by optimizing $\mf^2/\vf$. Such an optimization problem was already considered in \cite{Perry2018} for the spiked Wigner model. For the sake of completeness, we solve this problem by using the calculus of variations. Recall the density of the random variable $\sqrt{N}X_{ij}$ is $g$. 

To optimize $\f$, we need to maximize
\beq\label{target}
\left(\int_{-\infty}^{\infty}\f'(x)g(x)\dd x\right)^2\slash\left(\int_{-\infty}^{\infty}\f(x)^2g(x)\dd x\right)=\left(\int_{-\infty}^{\infty}\f(x)g'(x)\dd x\right)^2\slash\left(\int_{-\infty}^{\infty}\f(x)^2g(x)\dd x\right).
\eeq

Putting $(\f+\varepsilon\eta)$ in place of $\f$ in \eqref{target} and differentiating with respect to $\varepsilon$, we find that the optimal $\f$ satisfies
\beq\label{target2}
\left(\int_{-\infty}^{\infty}\eta(x)g'(x)\dd x\right)\left(\int_{-\infty}^{\infty}\f(x)^2g(x)\dd x\right)=\left(\int_{-\infty}^{\infty}\f(x)\eta(x)g(x)\dd x\right)\left(\int_{-\infty}^{\infty}\f(x)g'(x)\dd x\right)
\eeq
for any $\eta$. It is then easy to check that $\f=-Cg'/g$ is the only solution of \eqref{target2}. Since the value in \eqref{target} does not change if we replace $\f$ by $C\f$, and the effective SNR is increased with the entrywise transform $-g'/g$ is the optimal entrywise transformation for PCA.

\subsubsection{Multiplicative model} \label{subsec:variation_cov}
As we can see from the proof of Theorem \ref{thm:trans-cov} in Appendix \ref{sec:proof_trans-cov}, we need to maximize
\beq\begin{split}
	&\frac{2\left(\int_{-\infty}^{\infty}x\f(x)g(x)\dd x\right)\left(\int_{-\infty}^{\infty}\f'(x)g(x)\dd x\right)+\gamma\left(\int_{-\infty}^{\infty}\f'(x)g(x)\dd x\right)^2}{\left(\int_{-\infty}^{\infty}\f(x)^2g(x)\dd x\right)}
	\\&=\frac{-2\left(\int_{-\infty}^{\infty}x\f(x)g(x)\dd x\right)\left(\int_{-\infty}^{\infty}\f(x)g'(x)\dd x\right)+\gamma\left(\int_{-\infty}^{\infty}\f(x)g'(x)\dd x\right)^2}{\left(\int_{-\infty}^{\infty}\f(x)^2g(x)\dd x\right)}.
\end{split}\eeq
Putting $(\f+\varepsilon\eta)$ in place of $\f$ in \eqref{target} and differentiating with respect to $\varepsilon$, we find that the optimal $\f$ satisfies
\beq\label{opt}\begin{split}
	&-2\left(\int x\eta g\right)\left(\int \f g'\right)\left(\int \f^2 g\right)-2\left(\int x\f g\right)\left(\int \eta g'\right)\left(\int \f^2 g\right)+2\gamma\left(\int \eta g'\right)\left(\int \f g'\right)\left(\int \f^2 g\right)
	\\&+4\left(\int \f\eta g\right)\left(\int x\f g\right)\left(\int \f g'\right)-2\gamma\left(\int \f g'\right)^2\left(\int \f\eta g\right)=0
\end{split}\eeq
which is written with slight abuse of notation such as $\int x\eta g = \int_{-\infty}^{\infty}x \eta(x)g(x)\dd x$. Since the equation contains the terms
\[
\int x\eta g, \quad \int \eta g', \quad \int \f\eta g,
\]
it is natural to consider an ansatz
\beq
\f(x)=-\frac{g'(x)}{g(x)}+\alpha x
\eeq
for a constant $\alpha$. Collecting the terms involving $\int x\eta g$ and the terms involving $\int \eta g'$, we get
\[
2(\fg+\alpha)(\fg+2\alpha+\alpha^2)-4\alpha(1+\alpha)(\fg+\alpha)-2\alpha\gamma(\fg+\alpha)^2=0
\]
and
\[
-2(1+\alpha)(\fg+2\alpha+\alpha^2)-2\gamma(\fg+\alpha)(\fg+2\alpha+\alpha^2)+4(1+\alpha)(\fg+\alpha)+2\gamma(\fg+\alpha)^2=0.
\]
We can then check that 
\[
\alpha=\alpha_g=\frac{-\gamma \fg+\sqrt{4\fg+4\gamma \fg+\gamma^2\fg^2}}{2(1+\gamma)},
\]
and hence \eqref{opt} is satisfied with 
\[
\f(x)=-\frac{g'(x)}{g(x)}+\frac{-\gamma \fg+\sqrt{4\fg+4\gamma \fg+\gamma^2\fg^2}}{2(1+\gamma)} x.
\]
The corresponding effective SNR
\[
\lambda_{h_{\alpha_g}} \equiv \lambda_g= \gamma + \frac{\gamma^2 \fg}{2} + \frac{\gamma \sqrt{4\fg + 4\gamma \fg + \gamma^2 \fg^2}}{2}.
\] 

For a general $\alpha$, when the entrywise transform $h_{\alpha}$ is applied, the effective SNR
\[
\lambda_{h_{\alpha}} = \frac{2\gamma (1+\alpha)(\fg+\alpha) + \gamma^2 (\alpha+\fg)^2}{\alpha^2 + 2\alpha + \fg},
\]
In particular, if $\alpha = \sqrt{\fg}$,
\[
\lambda_{h_{\sqrt{\fg}}} = \gamma(1+\sqrt{\fg}) + \frac{\gamma^2}{2} (\fg+\sqrt{\fg}) \geq 2\gamma + \gamma^2 = \lambda
\]
where the inequality is strict if $\fg > 1$.

\subsection{Proof of Lemma \ref{lem:local_MP2}} \label{subsec:proof_lem}

\subsubsection{Proof of Lemma \ref{lem:local_MP2}}
Let $K=\F\caG X^T$ and $\caK=X\caG X^T$. Recall that $\sigma(z)=zs(z)+1$.
The key estimates in the proof of Lemma \ref{lem:local_MP2} are the following bounds on the entries of $K$ and $\caK$.
\begin{lem}\label{lem:entrywise}
	For $z \in \R$ outside an open interval containing $[d_-, d_+]$,
	\beq
	|K_{ij}-\widetilde{s}(z)\delta_{ij}|=\caO_\prec(N^{-1/2}), \qquad |\caK_{ij}-\check{s}(z)\delta_{ij}|=\caO_\prec(N^{-1/2}),
	\eeq
	where 
	\beq
	\widetilde{s}(z):=\sigma(z)\frac{\ef}{\sqrt{\vf}}, \qquad \check{s}(z):=zs(z) \left(\rat s(z)+\frac{\rat-1}{z}\right)^2 \frac{\ef^2}{\vf}+ \left(\rat s(z)+\frac{\rat-1}{z}\right).
	\eeq
\end{lem}
We postpone the proof of Lemma \ref{lem:entrywise} to Section \ref{subsec:proof_entrywise}. 

By definition,
\beq \label{eq:expand_K}
\langle\bsu,K\bsu\rangle- \widetilde{s}(z) \|\bsu \|^2=\sum_{i}u_i^2 (K_{ii}-\widetilde{s}(z))+\sum_{i\neq j}u_i K_{ij} u_j.
\eeq

From Lemma \ref{lem:entrywise} and the bound on $\| \bsu \|$, we find that the first term in the right side of \eqref{eq:expand_K} is $\caO_\prec(N^{-\phi})$. Applying Lemma \ref{lemma: LDE} with $q=N^{\phi}$ to the second term, 
\beq
\sum_{i\neq j}u_i K_{ij}u_j\prec \left(\frac1{M^2}\sum_{i\neq j}|K_{ij}|^2\right)^{1/2}+\frac{\max_{i\neq j}|K_{ij}|}{N^\phi}\prec N^{-\phi}.
\eeq
Applying Lemma \ref{lem:entrywise} again, we conclude that
\[
\langle\bsu,K\bsu\rangle = \widetilde{s}(z) \|\bsu \|^2 + \caO_\prec(N^{-\phi}) = \widetilde{s}(z) + \caO_\prec(N^{-\phi}).
\]
Since $K = \F\caG X^T$ is real symmetric, this proves the first part of the lemma. The second part of the lemma can be proved in the same manner with the estimates on the entries of $\caK$.

\subsubsection{Linearization} \label{subsec:linearization}

Recall that we have defined
\[
G \equiv G(z) = (\F\F^T -zI)^{-1}, \qquad \caG \equiv \caG(z) = (\F^T \F -zI)^{-1}.
\]
In the proof of Lemma \ref{lem:entrywise}, we use the formalism known as the linearization to simplify the computation. 
We define an $(M+N) \times (M+N)$ matrix $H_\F$ by
\beq \label{eq:linearization}
H_\F \equiv H_\F(z) =
\begin{pmatrix}
	-zI_M & \F \\
	\F^T & -I_N
\end{pmatrix},
\eeq
where $I_M$ and $I_N$ are the identity matrices with size $M$ and $N$, respectively. 

Let $R_\F(z) = H_\F(z)^{-1}$. (For the invertibility of $H_\F(z)$, we refer to Section 5.1 in \cite{Lee-Schnelli_Sample}.) By Schur's complement formula,
\beq \label{eq:resolvent_linearization}
R_\F(z) =
\begin{pmatrix}
	G(z) & G(z)\F \\
	\F^TG(z) & z\caG(z)
\end{pmatrix}.
\eeq
Therefore, 
\beq \label{eq:schur}
R_{ab}(z) = (\F\F^T -zI)^{-1}_{ab} = G_{ab}(z), \qquad R_{\alpha\beta}(z) = z(\F^T \F -zI)^{-1}_{\alpha-M, \beta-M} = z\caG_{\alpha-M, \beta-M}(z),
\eeq
and
\beq 
R_{\alpha a}(z)=R_{a\alpha}(z) = (G\F)_{a,\alpha-M}(z),
\eeq
where we use lowercase Latin letters $a, b, c, \dots$ for indices from $1$ to $M$ and Greek letters $\alpha, \beta, \gamma, \dots$ for indices from $(M+1)$ to $(M+N)$. We also use uppercase Latin letters $A, B, C, \dots$ for indices from $1$ to $(M+N)$. In the rest of Appendix \ref{sec:proof_trans}, we omit the subscript $Q$ for brevity.

For $\T \subset \{1, 2, \dots, M+N \}$, we define the matrix minor $H^{(\T)}$ by
\beq
(H^{(\T)})_{AB} := \mathbf{1}_{\{A, B \notin \T \}} H_{AB}\,.
\eeq
Moreover, for $A, B \notin \T$ we define
\beq
R^{(\T)}_{AB}(z) := (H^{(\T)})^{-1}_{AB},
\eeq
In the definitions above, we abbreviate $(\{A \})$ by $(A)$; similarly, we write $(AB)$ instead of $(\{A,B\})$.

We have the following identities for the matrix entries of $R$ and $R^{(\T)}$, which are elementary consequences of Schur's complement formula; see e.g.\ Lemma 5.1 of \cite{Lee-Schnelli_Sample}. 

\begin{lem}[Resolvent identities for $R$] \label{lem:res identity}
	Suppose that $z \in \R$ is outside an open interval containing $[d_-, d_+]$.
	
	- For $a \neq b$,
	\[
	R_{ab} = -R_{aa} \sum_{\alpha} H_{a \alpha} R^{(a)}_{\alpha b} = -R_{bb} \sum_{\beta} R^{(b)}_{a\beta} H_{\beta b}.
	\]
	
	- For $\alpha \neq \beta$,
	\[
	R_{\alpha\beta} = -R_{\alpha\alpha} \sum_{a} H_{\alpha a} R^{(\alpha)}_{a \beta} = -R_{\beta\beta} \sum_{b} R^{(\beta)}_{\alpha b} H_{b \beta}.
	\]
	
	- For any $a$ and $\alpha$,
	\[
	R_{a \alpha} = -R_{aa} \sum_{\beta} H_{a\beta} R^{(a)}_{\beta \alpha} = -R_{\alpha\alpha} \sum_{b} R^{(\alpha)}_{ab} H_{b \alpha}.
	\]
	
	- For $A, B \neq C$,
	\[
	R_{AB} = R_{AB}^{(C)} + \frac{R_{AC} R_{CB}}{R_{CC}}.
	\]
\end{lem}

Throughout this section, we will frequently use the estimate that all entries of $X$ and $\F$ (and hence all off-diagonal entries of $H$) are $\caO_{\prec}(N^{-1/2})$, which holds since all moments of the entries of $\sqrt{N} \F$ and $\sqrt{N} X$ are bounded. For the entries of $R$, we have the following estimates:

\begin{lem}\label{lem:entrywise_MPlaw}
	Let
	\beq
	\mathfrak{s}(z)=\left(\rat s(z)+\frac{\rat-1}{z}\right).
	\eeq
	For $z \in \R$ outside an open interval containing $[d_-, d_+]$,
	\beq
	\left|R_{ij}(z)-s(z)\delta_{ij}\right|, \left|R_{\mu\nu}(z)- z\mathfrak{s}(z)\delta_{\mu\nu}\right|, \left|R_{i \mu}(z)\right| \prec N^{-1/2}.
	\eeq
\end{lem}

\begin{proof}[Proof of Lemma \ref{lem:entrywise_MPlaw}]
	The first two estimates can be checked from Theorem 2.5 (and Remark 2.7) in \cite{Bloemendal-Erdos-Knowles-Yau-Yin2014} with the deterministic unit vectors $\mathbf{v}=\bse_i$ and $\mathbf{w}=\bse_j$ where $\bse_i\in\R^N$ or $\R^M$ is a standard basis vector whose $i$-th coordinate is 1 and all other coordinates are zero. For the last estimate, we apply Lemma \ref{lem:res identity} to find that
	\[
	R_{i \mu}(z) = -R_{ii} \sum_{\alpha} H_{i\alpha} R^{(i)}_{\alpha \mu}.
	\]
	Since $H_{i\alpha}$ and $R^{(i)}_{\alpha \mu}$ are independent, $R^{(i)}_{\alpha \mu} \prec N^{-1/2}$ for $\alpha \neq \mu$, and $R^{(i)}_{\mu\mu} = \Theta(1)$ with overwhelming probability, we find from Lemma \ref{lemma: LDE} that
	\[
	\sum_{\alpha} H_{i\alpha} R^{(i)}_{\alpha \mu} \prec \left( \frac{1}{N} \sum_{\alpha} | R^{(i)}_{\alpha \mu}|^2 \right)^{1/2} \prec N^{-1/2}.
	\]
\end{proof}

\subsubsection{Proof of Lemma \ref{lem:entrywise}} \label{subsec:proof_entrywise}

Throughout this section, for the sake of brevity, we will use the notation
\[
\F_{a\alpha} := \F_{a, (\alpha-M)} = H_{a\alpha}, \qquad X_{a\alpha} := X_{a, (\alpha-M)}.
\]

We begin by estimating the diagonal entry $K_{ii}$.
From Schur's complement formula, \eqref{eq:schur}, we can decompose it into
\beq \label{eq:decomp_K}
K_{ii} = (\F\caG X^T)_{ii} = \frac{1}{z} \sum_{\alpha} H_{i\alpha} R_{\alpha\alpha} X_{i\alpha} + \frac{1}{z} \sum_{\alpha\neq\beta} H_{i\alpha} R_{\alpha\beta} X_{i\beta}.
\eeq
From concentration inequalities it is not hard to see that
\[
\sum_{\alpha} \F_{i\alpha} X_{i\alpha} = \E[ \F_{i\alpha} X_{i\alpha} ] +\caO_{\prec}(N^{-1/2}) = \frac{\ef}{\sqrt{\vf}}+\caO_{\prec}(N^{-1/2}).
\]
Applying Lemma \ref{lem:entrywise_MPlaw}, we find for the first term in the right side of \eqref{eq:decomp_K} that
\beq \label{eq:decomp_K1}
\frac{1}{z} \sum_{\alpha} H_{i\alpha} R_{\alpha\alpha} X_{i\alpha} =\mathfrak{s}(z)\frac{\ef}{\sqrt{\vf}}+\caO_{\prec}(N^{-1/2}).
\eeq

We next estimate the second term in the right side of \eqref{eq:decomp_K}. We expand it with the resolvent identities in Lemma \ref{lem:res identity} as follows:
\beq \begin{split} \label{eq:decomp1}
	\sum_{\alpha\neq\beta} H_{i\alpha} R_{\alpha\beta} X_{i\beta} &= \sum_{\alpha\neq\beta} H_{i\alpha} R^{(i)}_{\alpha\beta} X_{i\beta} + \sum_{\alpha\neq\beta} H_{i\alpha} \frac{R_{\alpha i} R_{i\beta}}{R_{ii}} X_{i\beta} \\
	&= \sum_{\alpha\neq\beta} H_{i\alpha} R^{(i)}_{\alpha\beta} X_{i\beta} + \sum_{\alpha\neq\beta} H_{i\alpha} \frac{R_{\alpha i} R_{i\beta}}{s(z)} X_{i\beta} + \caO_{\prec}(N^{-1/2}).
\end{split} \eeq
Here, in the estimate for the second term, we simply counted the power (of $N$) as it involves two indices for the sum (hence $O(N^2)$ terms) of $H_{i\alpha}, R_{\alpha i}, R_{i\beta}, X_{i\beta} \prec N^{-1/2}$, hence $\sum_{\alpha\neq\beta} H_{i\alpha} R_{\alpha i} R_{i\beta} X_{i\beta} = \caO_{\prec}(1)$.
Applying Lemma \ref{lemma: LDE} to the first term in the right side of \eqref{eq:decomp1},
\[
\sum_{\alpha\neq\beta} H_{i\alpha} R^{(i)}_{\alpha\beta} X_{i\beta} \prec \left( \frac{1}{N^2} \sum_{\alpha, \beta} |R^{(i)}_{\alpha\beta}|^2 \right)^{1/2} \prec N^{-1/2}.
\]
For the second term in the right side of \eqref{eq:decomp1}, we further expand it to find
\[ \begin{split}
\sum_{\alpha\neq\beta} H_{i\alpha} R_{\alpha i} R_{i\beta} X_{i\beta} = \sum_{\alpha\neq\beta} H_{i\alpha} R_{\alpha i} R_{i\beta} X_{i\beta} = -\sum_{\alpha\neq\beta} H_{i\alpha} \left( R_{ii} \sum_{\mu} R^{(i)}_{\alpha\mu} H_{\mu i} R_{i\beta} X_{i\beta} \right)
\end{split} \]
Note that
\[
\sum_{\mu} R^{(i)}_{\alpha\mu} H_{\mu i} \prec N^{-1/2},
\]
as in the proof of Lemma \ref{lem:entrywise_MPlaw}. Since 
\[
|R_{ij} - s(z)| \prec N^{-1/2}, \quad R_{i\beta} = R^{(\alpha)}_{i\beta} + \frac{R_{i\alpha} R_{\alpha\beta}}{R_{\alpha\alpha}} = R^{(\alpha)}_{i\beta} + N^{-1},
\]
we have
\beq \label{eq:expand_oneside} \begin{split}
	&-\sum_{\alpha\neq\beta} H_{i\alpha} \left( R_{ii} \sum_{\mu} R^{(i)}_{\alpha\mu} H_{\mu i} R_{i\beta} X_{i\beta} \right) = -s(z) \sum_{\alpha\neq\beta} H_{i\alpha} \left( \sum_{\mu} R^{(i)}_{\alpha\mu} H_{\mu i} R^{(\alpha)}_{i\beta} X_{i\beta} \right) + \caO_{\prec} (N^{-1/2}) \\
	&= -s(z) \sum_{\alpha\neq\beta} H_{i\alpha} \left( \sum_{\mu:\mu \neq \alpha} R^{(i)}_{\alpha\mu} H_{\mu i} R^{(\alpha)}_{i\beta} X_{i\beta} \right) -s(z) \sum_{\alpha\neq\beta} (H_{i\alpha})^2 R^{(i)}_{\alpha\alpha} R^{(\alpha)}_{i\beta} X_{i\beta} + \caO_{\prec} (N^{-1/2}).
\end{split} \eeq
Applying Lemma \ref{lemma: LDE} again to the first term in the right side of \eqref{eq:expand_oneside},
\[ \begin{split}
\sum_{\alpha\neq\beta} H_{i\alpha} \left( \sum_{\mu:\mu \neq \alpha} R^{(i)}_{\alpha\mu} H_{\mu i} R^{(\alpha)}_{i\beta} X_{i\beta} \right) &\prec \left( \frac{1}{N} \sum_{\alpha} \left| \sum_{\beta: \beta\neq \alpha} \left[ \sum_{\mu:\mu \neq \alpha} R^{(i)}_{\alpha\mu} H_{\mu i} \right] R^{(\alpha)}_{i\beta} X_{i\beta} \right|^2 \right)^{1/2} \\
& \prec \left( \frac{1}{N} \sum_{\alpha} \left[ \sum_{\beta: \beta\neq \alpha} N^{-1/2} \left| R^{(\alpha)}_{i\beta} X_{i\beta} \right| \right]^2 \right)^{1/2} \prec N^{-1/2}.
\end{split} \]
Similarly, by expanding $R^{(\alpha)}_{i\beta}$, we find for the second term in the right side of \eqref{eq:expand_oneside} that
\[ \begin{split}
-&s(z) \sum_{\alpha\neq\beta} (H_{i\alpha})^2 R^{(i)}_{\alpha\alpha} R^{(\alpha)}_{i\beta} X_{i\beta} = zs(z) \mathfrak{s}(z) \sum_{\alpha\neq\beta} (H_{i\alpha})^2 R^{(\alpha)}_{ii} \sum_{\nu:\nu \neq \alpha}H^{(\alpha)}_{i\nu} R^{(i\alpha)}_{\nu\beta} X_{i\beta} + \caO_{\prec} (N^{-1/2}) \\
&= zs(z)^2 \mathfrak{s}(z) \sum_{\alpha\neq\beta} (H_{i\alpha})^2 \sum_{\nu:\nu \neq \alpha, \beta} H_{i\nu} R^{(i\alpha)}_{\nu\beta} X_{i\beta} + zs(z)^2 \mathfrak{s}(z) \sum_{\alpha\neq\beta} (H_{i\alpha})^2 H_{i\beta} R^{(i\alpha)}_{\beta\beta} X_{i\beta} + \caO_{\prec} (N^{-1/2}) \\
&= z^2 s(z)^2 \mathfrak{s}(z)^2 \sum_{\alpha\neq\beta} (H_{i\alpha})^2 H_{i\beta} X_{i\beta} + \caO_{\prec} (N^{-1/2}),
\end{split} \]
where we used Lemma \ref{lemma: LDE} to find
\[
\sum_{\nu \neq \beta:\nu, \beta \neq \alpha} H_{i\nu} R^{(i\alpha)}_{\nu\beta} X_{i\beta} \prec \left( \frac{1}{N^2} \sum_{\nu \neq \beta:\nu, \beta \neq \alpha} \left| R^{(i\alpha)}_{\nu\beta} \right|^2 \right)^{1/2} \prec N^{-1/2}.
\]
Thus,
\[
\sum_{\alpha\neq\beta} H_{i\alpha} R_{\alpha i} R_{i\beta} X_{i\beta} = z^2 s(z)^2 \mathfrak{s}(z)^2 \sum_{\alpha\neq\beta} (H_{i\alpha})^2 H_{i\beta} X_{i\beta} + \caO_{\prec} (N^{-1/2}) = z^2 s(z)^2 \mathfrak{s}(z)^2 \frac{\ef}{\sqrt{\vf}} + \caO_{\prec} (N^{-1/2}),
\]
and putting it back to \eqref{eq:decomp1} and \eqref{eq:decomp_K}, together with \eqref{eq:decomp_K1}, we conclude that
\beq \label{eq:estimate_K}
K_{ii} = \mathfrak{s}(z) \frac{\ef}{\sqrt{\vf}} + z s(z) \mathfrak{s}(z)^2 \frac{\ef}{\sqrt{\vf}} + \caO_{\prec} (N^{-1/2}) = \sigma(z) \frac{\ef}{\sqrt{\vf}} + \caO_{\prec} (N^{-1/2}),
\eeq
where we used the identity $zs(z)\mathfrak{s}(z) = -\sigma(z)$.
In the same manner, we also find that
\beq \begin{split} \label{eq:estimate_caK}
	\caK_{ii} &= \frac{1}{z} \sum_{\alpha} X_{i\alpha} R_{\alpha\alpha} X_{i\alpha} + zs(z)\mathfrak{s}(z)^2 \sum_{\alpha\neq\beta} X_{i\alpha} H_{i\alpha} H_{i\beta} X_{i\beta} + \caO_{\prec} (N^{-1/2}) \\
	&= \mathfrak{s}(z) + zs(z)\mathfrak{s}(z)^2 \frac{\ef^2}{\vf} + \caO_{\prec} (N^{-1/2}).
\end{split} \eeq

We next estimate the off-diagonal entry $K_{ij}$. We expand it as
\beq \begin{split} \label{eq:decomp_off}
	K_{ij} &= \frac{1}{z} \sum_{\alpha, \beta} H_{i\alpha} R_{\alpha\beta} X_{j\beta} = \frac{1}{z} \sum_{\alpha, \beta} H_{i\alpha} R^{(i)}_{\alpha\beta} X_{j\beta} + \frac{1}{z} \sum_{\alpha, \beta} H_{i\alpha} \frac{R_{\alpha i} R_{i\beta}}{R_{ii}} X_{j\beta} \\
	&= \frac{1}{z} \sum_{\alpha, \beta} H_{i\alpha} R^{(ij)}_{\alpha\beta} X_{j\beta} + \frac{1}{z} \sum_{\alpha, \beta} H_{i\alpha} \frac{R^{(i)}_{\alpha j} R^{(i)}_{j\beta}}{R^{(i)}_{jj}} X_{j\beta} + \frac{1}{z} \sum_{\alpha, \beta} H_{i\alpha} \frac{R^{(j)}_{\alpha i} R^{(j)}_{i\beta}}{R^{(i)}_{jj}} X_{j\beta} + \caO_{\prec}(N^{-1/2})
\end{split} \eeq
From Lemma \ref{lemma: LDE},
\[
\sum_{\alpha, \beta} H_{i\alpha} R^{(ij)}_{\alpha\beta} X_{j\beta} \prec N^{-1/2}.
\]
We also have
\[
\sum_{\alpha, \beta} H_{i\alpha} \frac{R^{(i)}_{\alpha j} R^{(i)}_{j\beta}}{R^{(i)}_{jj}} X_{j\beta} \prec \left( \frac{1}{N} \sum_{\alpha} \left| \sum_{\beta} \frac{R^{(i)}_{\alpha j} R^{(i)}_{j\beta}}{R^{(i)}_{jj}} X_{j\beta} \right|^2 \right)^{1/2} \prec \left( \frac{1}{N} \sum_{\alpha} \left| \sum_{\beta} N^{-3/2} \right|^2 \right)^{1/2} \prec N^{-1/2}
\]
and a similar estimate holds for the third term in the right side of \eqref{eq:decomp_off}. Thus,
\[
K_{ij} \prec N^{-1/2}
\]
In the same manner, we also find that $\caK_{ij} \prec N^{-1/2}$. Together with \eqref{eq:estimate_K} and \eqref{eq:estimate_caK}, this proves Lemma \ref{lem:entrywise}.

\section{Proof of Theorem \ref{thm:CLT}} \label{sec:proof_CLT}

Recall that for a function $f$ analytic on an open set containing an interval $[d_-, d_+]$
\beq 
\sum_{i=1}^M f(\mu_i) = \sum_{i=1}^M \frac{1}{2\pi\ii} \oint_{\Gamma} \frac{f(z)}{z-\mu_i} \dd z = -\frac{1}{2\pi\ii} \oint_{\Gamma} f(z) \Tr (YY^T -zI)^{-1} \dd z
\eeq
for any contour $\Gamma$ containing $\mu_1, \mu_2, \dots, \mu_N$. Our goal is to track the change of the LSS by finding the change of the trace of the resolvent $\Tr (YY^T -zI)^{-1}$ and conclude that the change is decomposed into the deterministic part and the random part, where the latter converges to $0$ with overwhelming probability. We also directly compute the change of the mean from the null model to the non-null model.

\subsection{Additive model}\label{subsec:spikedmean}
Let
\beq 
Y(\theta)=\theta\sqrt{\lambda}\bsu\bsv^T+X
\eeq
for $\theta \in [0, 1]$. Note that $Y(0) = X$ and $Y(1) = Y$. Denote by $\mu_1(\theta) \geq \mu_2(\theta) \geq \dots \geq \mu_M(\theta)$ the eigenvalues of $Y(\theta) Y(\theta)^T$. We also define the resolvent
\beq
G(\theta, z) = (Y(\theta)Y(\theta)^T -zI)^{-1}, \qquad \caG(\theta, z) = (Y(\theta)^T Y(\theta) -zI)^{-1}
\eeq
for $z \in \C$.

We choose ($N$-independent) constants $a_- <d_-$, $a_+ > d_+$, and $v_0 \in (0, 1)$ so that the function $f$ is analytic on the rectangular contour $\Gamma$ whose vertices are $(a_- \pm \ii v_0)$ and $(a_+ \pm \ii v_0)$. With overwhelming probability, all eigenvalues of $Y(\theta)Y(\theta)^T$ are contained in $\Gamma$. Applying Cauchy's integral formula, as in \eqref{eq:Cauchy}, we find that
\beq
\sum_{i=1}^M f(\mu_i(1)) - \sum_{i=1}^M f(\mu_i(0)) = -\left(\frac{1}{2\pi \ii} \oint_{\Gamma}f(z) \left( \Tr G(1,z) - \Tr G(0,z) \right) \dd z \right)
\eeq
To estimate the difference $\Tr G(1,z) - \Tr G(0,z)$, we consider its derivative $\frac{\partial}{\partial \theta}\Tr G(\theta,z)$. Note that
\beq
\frac{\partial G_{ab}(\theta)}{\partial Y_{jk}(\theta)}=-G_{aj}(\theta)(Y(\theta)^T G(\theta))_{kb}-(G(\theta)Y(\theta))_{ak}G_{jb}(\theta), \qquad 
\frac{\dd Y_{jk}(\theta)}{\dd\theta}=\sqrt{\lambda}u_j v_k.
\eeq
Thus, by chain rule
\beq \begin{split} \label{eq:der_Tr_mean}
	\frac{\partial}{\partial \theta}\Tr G(\theta,z)&=\sum_{a=1}^{M}\sum_{j=1}^{M}\sum_{k=1}^{N}\frac{\partial Y_{jk}(\theta)}{\partial\theta}\frac{\partial G_{aa}(\theta)}{\partial Y_{jk}(\theta)} \\
	&=-\sum_{a=1}^{M}\sum_{j=1}^{M}\sum_{k=1}^{N}\sqrt{\lambda}u_jv_k[G_{aj}(\theta)(Y(\theta)^T G(\theta))_{ka}+(G(\theta)Y(\theta))_{ak}G_{ja}(\theta)] \\
	&=-2\sum_{a=1}^{M}\sum_{j=1}^{M}\sum_{k=1}^{N}\sum_{\ell=1}^{M}\sqrt{\lambda}u_jv_k[Y_{\ell k}(\theta)G_{\ell a}(\theta)G_{aj}(\theta)]
\end{split} \eeq
From the fact 
\[
\left(\frac{\partial}{\partial z}G(\theta)\right)_{\ell j}=(G(\theta)^2)_{\ell j}=\sum_a G_{\ell a}(\theta)G_{aj}(\theta),
\]
we then find that
\beq \begin{split} \label{eq:Tr_derivative}
	\frac{\partial}{\partial \theta}\Tr G(\theta,z) =-2\sqrt{\lambda}\frac{\partial}{\partial z}\sum_{j=1}^{M}\sum_{k=1}^{N}u_jv_k(G(\theta)Y(\theta))_{jk}=-2\sqrt{\lambda}\frac{\partial}{\partial z}\langle \bsu,G(\theta)Y(\theta)\bsv\rangle.
\end{split} \eeq

It remains to estimate $\frac{\partial}{\partial z}\langle \bsu,G(\theta)Y(\theta)\bsv\rangle$. Note that
\[
\langle \bsu,G(\theta)Y(\theta)\bsv\rangle=\theta\sqrt{\lambda}\langle \bsu,G(\theta)\bsu\rangle+\langle \bsu,G(\theta)X\bsv\rangle.
\]
We consider the resolvent expansion 
\beq\label{eq:resolvent_mean}\begin{split}
	G(0,z)-G(\theta,z)&=G(\theta,z)\,(H(\theta)-H(0))\,G(0,z) \\
	&=G(\theta,z)\,(\theta^2\lambda \bsu\bsu^T+\theta\sqrt{\lambda}X\bsv\bsu^T+\theta\sqrt{\lambda}\bsu\bsv^TX^T)\,G(0,z).
\end{split}\eeq
Taking inner products with $\bsu$ and $\bsv$, we obtain
\beq\begin{split}\label{eq:ex1}
	\langle \bsu,G(0)\bsu\rangle&=\langle \bsu,G(\theta)\bsu\rangle+\theta^2\lambda\langle \bsu,G(\theta)\bsu\rangle\langle \bsu,G(0)\bsu\rangle\\&~~~+\theta\sqrt{\lambda}\langle \bsu,G(\theta)X\bsv\rangle\langle \bsu,G(0)\bsu\rangle+\theta\sqrt{\lambda}\langle\bsu,G(0)X\bsv\rangle\langle \bsu,G(\theta)\bsu\rangle
\end{split}\eeq
and
\beq\begin{split}\label{eq:ex2}
	\langle \bsu,G(0)X\bsv\rangle&=\langle \bsu,G(\theta)X\bsv\rangle+\theta^2\lambda\langle \bsu,G(\theta)X\bsv\rangle\langle \bsu,G(0)X\bsv\rangle\\&~~~+\theta\sqrt{\lambda}\langle \bsu,G(\theta)X\bsv\rangle\langle \bsu,G(0)X\bsv\rangle+\theta\sqrt{\lambda}\langle \bsv,X^T G(0)X\bsv\rangle\langle \bsu,G(\theta)\bsu\rangle,
\end{split}\eeq
where we omitted $z$-dependence for brevity. We then use the following result to control the terms in \eqref{eq:ex1} and \eqref{eq:ex2}. Recall the definition of $s(z)$ and $\mathfrak{s}(z)$ in Lemmas \ref{lem:local_MP} and \ref{lem:entrywise_MPlaw}. Moreover, we consider the same linearization $H_X(z)$ of the matrix $X$ and the corresponding resolvent $R_{X}(z)=H_{X}(z)^{-1}.$ as in \eqref{eq:linearization} and \eqref{eq:resolvent_linearization}.

\begin{lem}[Isotropic local law] \label{lem:local_law}
	For an $N$-independent constant $\varepsilon > 0$, let $\Gamma^{\varepsilon}$ be the $\varepsilon$-neighborhood of $\Gamma$, i.e.,
	\[
	\Gamma^{\varepsilon} = \{ z \in \mathbb{C} : \min_{w \in \Gamma} |z-w| \leq \varepsilon \}.
	\]
	Choose $\varepsilon$ small so that the distance between $\Gamma^{\varepsilon}$ and $[d_-, d_+]$ is larger than $2\varepsilon$, i.e., 
	\beq
	\min_{w \in \Gamma^{\varepsilon}, x \in [d_-, d_+]} |x-w| > 2\varepsilon.
	\eeq
	Then, for any unit vectors $\bsx,\bsy\in\mathbb{C}^{M+N}$ independent of $X$,
	\beq \begin{split} \label{eq:iso_spike}
		\left|\left\langle \bsx,(R_X(z)-\Pi(z)) \bsy\right\rangle\right|\prec N^{-1/2}, 
	\end{split} \eeq
	uniformly on $z\in \Gamma^\varepsilon$, where
	\beq 
	\Pi(z) =
	\begin{pmatrix}
		s(z)\cdot I_M & 0 \\
		0 & z\frak{s}(z)\cdot I_N
	\end{pmatrix}.
	\eeq
\end{lem}
\begin{proof}
	See Theorems 3.6, 3.7, Corollary 3.9, and Remark 3.10 in \cite{Knowles-Yin2016}. Note that $\im \mathfrak{s}(z), \im s(z) = \Theta(\eta)$ on the vertical part of $\Gamma_{\varepsilon}$, i.e., the neighborhood of the line segment joining $(a_+ + \ii v_0)$ and $(a_+ - \ii v_0)$ (respectively $(a_- + \ii v_0)$ and $(a_- - \ii v_0)$).
\end{proof}

Set
\[\begin{split}
A:=\langle \bsu,G(0, z)\bsu\rangle,\qquad
B:=\langle \bsu,G(0, z)X\bsv\rangle,\qquad
C:=\langle \bsv,X^T G(0, z)X\bsv\rangle.
\end{split}\]
Recall that 
\beq 
R_X(z) =
\begin{pmatrix}
	G(0,z) & G(0,z)X \\
	X^TG(0,z) & z\caG(0,z)
\end{pmatrix}.
\eeq
Then, as consequences of Lemma \ref{lem:local_law} with appropriate choices of the deterministic vectors,
\beq
A= s(z) +\caO_\prec(N^{-1/2}), \qquad C =\langle \bsv, z\caG(0, z)\bsv\rangle+1+\caO(N^{-1/2})=\rat(zs(z)+1)+\caO_\prec(N^{-1/2}),
\eeq
and
\[
B=\caO_\prec(N^{-1/2}).
\]

We thus have from \eqref{eq:ex1} and \eqref{eq:ex2} that
\beq \begin{split} \label{eq:iso_perturb}
	\langle \bsu,G(\theta)X\bsv\rangle &= -\frac{\theta\rat \sqrt{\lambda}s(z)(zs(z)+1)}{\theta^2\lambda zs(z)+\theta^2\lambda+1}+\caO_\prec(N^{-1/2}) \\
	\langle \bsu,G(\theta)\bsu\rangle &=\frac{s(z)}{\theta^2\lambda zs(z)+\theta^2\lambda+1}+\caO_\prec(N^{-1/2})
\end{split} \eeq
and hence
\beq
\langle \bsu,G(\theta)Y(\theta)\bsv\rangle =\theta\sqrt{\lambda}\langle \bsu,G(\theta)\bsu\rangle+\langle \bsu,G(\theta)X\bsv\rangle =\frac{\theta\sqrt{\lambda} zs(z)+\theta\sqrt{\lambda}}{\theta^2\lambda zs(z)+\theta^2\lambda+1}+\caO_\prec(N^{-1/2}).
\eeq
Note that this estimate is uniform on $\theta$.
Differentiating it with respect to $z$ and plugging it back to \eqref{eq:Tr_derivative}, we get
\[
\frac{\partial}{\partial \theta}\Tr G(\theta,z)=-\frac{2\theta\lambda\frac{\dd}{\dd z}(zs(z)+1)}{(\theta^2\lambda zs(z) +\theta^2\lambda+1)^2}+\caO_\prec(N^{-1/2})
\] 
and, integrating over $\theta$, we obtain
\beq \label{eq:difference_Tr}
\Tr G(1,z)-\Tr G(0,z)=\int_{0}^{1}\frac{\partial}{\partial \theta}\Tr G(\theta,z) \dd\theta=-\frac{\frac{\dd}{\dd z}\lambda(zs(z)+1)}{\lambda zs(z)+\lambda+1}+\caO_\prec(N^{-1/2}).
\eeq

We now invoke the following relation between the Marchenko--Pastur law and the Wigner semicircle law. Let 
\[
s_{sc}(z) = \frac{-z+\sqrt{z^2 -4}}{2}
\]
be the Stieltjes transform of the Wigner semicircle law and 
\[
\varphi(z)=\frac1{\sqrt{\rat}}(z-(1+\rat)).
\]
Then
\[
\sqrt{\rat} (zs(z)+1)=s_{sc}(\varphi(z)).
\]
We thus have
\beq \begin{split}
	\frac1{2\pi \ii}\oint_\Gamma f(z)\frac{\lambda\frac{\dd}{\dd z}(zs(z)+1)}{\lambda zs(z)+\lambda+1} \dd z &=\frac1{2\pi \ii}\oint_\Gamma (f\circ\phi)(\varphi(z))\frac{\lambda s_{sc}'(\varphi(z))\varphi'(z)}{\lambda s_{sc}(\varphi(z))+\sqrt{\rat}} \dd z \\
	&=\frac1{2\pi \ii}\oint_{\wt\Gamma} \wt{f}(\varphi)\frac{\lambda s_{sc}'(\varphi)}{\lambda s_{sc}(\varphi)+\sqrt{\rat}} \dd\varphi
\end{split} \eeq
where we let $(f\circ\phi)(z)=\wt{f}(z)$ and $\wt\Gamma=\varphi(\Gamma)$. (Note that $\wt \Gamma$ contains the interval $[-2,2]$.)

So far, we have proved that
\beq \label{eq:mean_mean}
\sum_{i=1}^M f(\mu_i(1)) - \sum_{i=1}^M f(\mu_i(0)) = \frac1{2\pi \ii}\oint_{\wt\Gamma} \wt{f}(\varphi)\frac{\lambda s_{sc}'(\varphi)}{\lambda s_{sc}(\varphi)+\sqrt{\rat}} \dd\varphi+\caO_\prec(N^{-1/2}). 
\eeq
Since the difference in \eqref{eq:mean_mean} is the sum of a deterministic term and a random term stochastically dominated by $N^{-1/2}$, we can see that the CLT holds for the LSS with the non-null model $Y(1)$. Moreover, the variance is the same as that of the null model, which is
\beq
V_Y(f) =2\sum_{\ell=1}^\infty \ell\tau_\ell(\wt{f})^2+(w_4-3)\tau_1(\wt{f})^2.
\eeq
(See, e.g., \cite{Baik-Lee2018}.)

The change of the mean is the first term in the right side of \eqref{eq:mean_mean}, which can be computed by following the proof of Lemma 4.4 in \cite{Baik-Lee2017}. We obtain
\beq
m_Y(f) = \frac{\wt{f}(2) + \wt{f}(-2)}{4}  -\frac{1}{2} \tau_0(\wt{f})-(w_4-3)\tau_2(\wt{f})+\sum_{\ell=1}^{\infty}\left(\frac{\lambda}{\sqrt{\rat}}\right)^{\ell}\tau_\ell(\wt{f}).
\eeq
This proves the first part of Theorem \ref{thm:CLT} for the additive model.

\subsection{Multiplicative model}
For the multiplicative model, we will follow the same strategy as in the additive model. Let
\beq 
Y(\theta)= X+ \theta\gamma \bsu\bsu^T X
\eeq
for $\theta \in [0, 1]$. Note that $Y(0) = X$ and $Y(1) = Y$. As in Section \ref{subsec:spikedmean}, we denote by $\mu_1(\theta) \geq \mu_2(\theta) \geq \dots \geq \mu_M(\theta)$ the eigenvalues of $Y(\theta) Y(\theta)^T$, and also let
\beq
G(\theta, z) = (Y(\theta)Y(\theta)^T -zI)^{-1}, \qquad \caG(\theta, z) = (Y(\theta)^T Y(\theta) -zI)^{-1}
\eeq
for $z \in \C$. We have the relations
\beq \begin{split}
	\frac{\partial G_{ab}(\theta)}{\partial Y_{jk}(\theta)}=-G_{aj}(\theta)(Y(\theta)^T G(\theta))_{kb}-(G(\theta)Y(\theta))_{ak}G_{jb}(\theta), \qquad 
	\frac{\partial Y_{jk}(\theta)}{\partial\theta}=\gamma\sum_{\ell=1}^{M}u_ju_\ell X_{\ell k}.
\end{split} \eeq

Following \eqref{eq:der_Tr_mean}-\eqref{eq:Tr_derivative}, we get
\beq \begin{split}
	\frac{\partial}{\partial \theta}\Tr G(\theta,z) &= -\gamma\sum_{a=1}^{M}\sum_{j=1}^{M}\sum_{k=1}^{N}\sum_{\ell=1}^{M}u_ju_\ell X_{\ell k}[G_{aj}(\theta)(Y(\theta)^T G(\theta))_{ka}+(G(\theta)Y(\theta))_{ak}G_{ja}(\theta)] \\
	&=-2\gamma\sum_{a=1}^{M}\sum_{j=1}^{M}\sum_{k=1}^{N}\sum_{s=1}^{M}\sum_{\ell=1}^{M}u_ju_\ell X_{\ell k}[Y_{s k}(\theta)G_{s a}(\theta)G_{aj}(\theta)] \\
	&=-2\gamma\frac{\partial}{\partial z}\sum_{j=1}^{M}\sum_{k=1}^{N}\sum_{\ell=1}^{M}u_ju_\ell X_{\ell k}(G(\theta)Y(\theta))_{jk} \\
	&=-2\gamma\frac{\partial}{\partial z}\langle \bsu,G(\theta)Y(\theta)X^T \bsu\rangle=-2\gamma\frac{\partial}{\partial z}\langle \bsu,G(\theta)Y(\theta)Y(0)^T \bsu\rangle.
\end{split} \eeq
Moreover, since 
\beq
Y(0) = X = (I + \theta\gamma \bsu\bsu^T)^{-1} Y(\theta) = \left( I-\frac{\theta\gamma}{1+\theta\gamma}\bsu\bsu^T \right) Y(\theta),
\eeq  
we have
\beq \begin{split}
	\langle \bsu,G(\theta)Y(\theta)Y(0)^T \bsu\rangle &=\langle \bsu,G(\theta)Y(\theta)Y(\theta)^T (I + \theta\gamma \bsu\bsu^T)^{-1} \bsu\rangle =\langle \bsu,(I+zG(\theta))(I + \theta\gamma \bsu\bsu^T)^{-1} \bsu\rangle \\
	&=\frac{1}{1+\theta\gamma}+\frac{z}{1+\theta\gamma}\langle \bsu,G(\theta)\bsu\rangle.\label{eq:term_cov}
\end{split} \eeq
To estimate the term $\langle \bsu,G(\theta)\bsu\rangle$, we use the following Anisotropic local law in \cite{Knowles-Yin2016}.

\begin{lem}[Anisotropic local law] \label{lem:anisotropic}
	Let $\Gamma^{\varepsilon}$ be the $\varepsilon$-neighborhood of $\Gamma$ as in Lemma \ref{lem:local_law}. Then, for any unit vectors  $\bsx,\, \bsy \in \C^M$ independent of $X$, the following estimate holds uniformly on $z \in \Gamma^{\varepsilon}$ :
	\beq \label{eq:aniso_spike}
	\left|\left\langle \wt \bsx, \left( G(\theta, z)+ \left( zI + z\mathfrak{s}(z) (I + \theta \gamma \bsu \bsu^T)^2 \right)^{-1} \right) \wt \bsy\right\rangle \right| \prec N^{-\frac{1}{2}}.
	\eeq
\end{lem}
\begin{proof}
	The proof of Lemma \ref{lem:anisotropic} is the same as that of Lemma \ref{lem:local_law}.
\end{proof}
From Lemma \ref{lem:anisotropic}, we find that
\beq
\langle \bsu,G(\theta)\bsu\rangle = -\left\langle \bsu, \left( zI + z\mathfrak{s}(z) (I + \theta \gamma \bsu \bsu^T)^2 \right)^{-1} \bsu \right\rangle  +\caO(N^{-1/2}) = -\frac{1}{(1+\theta\gamma)^2 z(1+\mathfrak{s}(z) )}+\caO(N^{-1/2}),
\eeq
and plugging it into \eqref{eq:term_cov}, we obtain
\beq\begin{split}
	\langle \bsu,G(\theta)Y(\theta)Y(0)^T \bsu\rangle=\frac{1}{1+\theta\gamma}-\frac{1}{(1+\theta\gamma)(1+ (1+\theta\gamma)^2 \mathfrak{s}(z) )}+\caO(N^{-1/2}).
\end{split}
\eeq
We thus get
\beq
\frac{\partial}{\partial \theta}\Tr G(\theta,z)=-2\gamma\frac{(1+\theta\gamma) \mathfrak{s}'(z)}{(1+ (1+\theta\gamma)^2 \mathfrak{s}(z) )^2}+\caO(N^{-1/2}),
\eeq
and integrating it yields 
\beq \label{eq:difference_Tr_2}
\Tr G(1,z)-\Tr G(0,z)=-\frac{\lambda\mathfrak{s}'(z)}{(1+\mathfrak{s}(z) )(1+(1+\lambda)\mathfrak{s}(z))}+\caO(N^{-1/2})=-\frac{\lambda \frac{\dd}{\dd z}(zs(z)+1)}{\lambda zs(z)+\lambda+1}+\caO(N^{-1/2}).
\eeq

Since \eqref{eq:difference_Tr_2} coincides with \eqref{eq:difference_Tr}, the rest of the proof is exactly the same as in the additive case. This finishes the proof of the first part of Theorem \ref{thm:CLT}.

\subsection{Computation of the test statistic} \label{sec:compute}
In this section, we prove the second part of Theorem \ref{thm:CLT} and also provide the details on the computation of the test statistic in Theorem \ref{thm:test}. Recall that
\beq
m_Y(f)|_{\bsH_1} - m_Y(f)|_{\bsH_0} = \sum_{\ell=1}^{\infty}\left(\frac{\SNR}{\sqrt{\rat}}\right)^{\ell}\tau_\ell(\wt{f})
\eeq
and
\beq \begin{split}
	V_{Y}(f) = 2\sum_{\ell=2}^\infty \ell\tau_\ell(\wt{f})^2+(w_4-1)\tau_1(\wt{f})^2.
\end{split} \eeq
Assuming $w_2 > 0$ and $w_4 > 1$, from Cauchy's inequality and the identity $\log(1-\lambda) = -\sum_{\ell=1}^{\infty} \lambda^{\ell}/\ell$,
\beq\begin{split} \label{eq:Cauchy_ineq}
	\left|\frac{m_Y(f)|_{\bsH_1} - m_Y(f)|_{\bsH_0}}{\sqrt{V_Y(f)}}\right|^2 &\leq  \frac{\SNR^2}{\rat(w_4-1)}-\frac{\SNR^2}{2\rat} + \frac1{2}\sum_{\ell=1}^{\infty} \frac{1}{\ell}\left(\frac{\SNR^2}{\rat}\right)^{\ell} 
	\\&= \frac{\SNR^2}{\rat}\left( \frac{1}{w_4-1} - \frac{1}{2} \right)  - \frac{1}{2}\displaystyle \log\left(1-\frac{\SNR^2}{\rat}\right) = \left|\frac{m(w)-m(0)}{\sqrt{V_0}}\right|^2,
\end{split}\eeq
which proves the first part of the theorem. The equality in \eqref{eq:Cauchy_ineq} holds if and only if
\beq \label{eq:Cauchy_equal}
\sqrt{\rat}(w_4-1)\tau_1(\wt{f}) = \frac{2\ell (\sqrt{\rat})^\ell\tau_\ell(\wt{f})}{\SNR^{\ell-1}} \qquad (\ell = 2, 3, 4, \dots).
\eeq

We now find all functions $f$ that satisfy \eqref{eq:Cauchy_equal}. Letting $2C$ be the common value in \eqref{eq:Cauchy_equal},
\beq \label{eq:tau_ell}
\tau_1(\wt{f}) = \frac{2C \SNR }{\sqrt{\rat}(w_4-1)}, \quad \tau_{\ell}(\wt{f}) = \frac{C\SNR^{\ell}}{\ell(\sqrt{\rat})^\ell} \qquad (\ell = 2, 3, 4, \dots).
\eeq
We can expand $\wt{f}$ in terms of the Chebyshev polynomials as
\beq
\wt f(x) = \sum_{\ell=0}^{\infty} C_{\ell} T_{\ell} \left( \frac{x}{2} \right).
\eeq
From the orthogonality relation of the Chebyshev polynomials, we get for $\ell \geq 1$ that
\beq
\tau_{\ell}(\wt f) = \frac{C_{\ell} }{\pi} \int_{-2}^2 T_{\ell} \left( \frac{x}{2} \right) T_{\ell} \left( \frac{x}{2} \right) \frac{\dd x}{\sqrt{4-x^2}} = \frac{C_{\ell} }{\pi} \int_{-1}^1 T_{\ell} \left( y \right) T_{\ell} \left( y \right) \frac{\dd y}{\sqrt{1-y^2}} = \frac{C_{\ell}}{2}.
\eeq
Thus, \eqref{eq:tau_ell} holds if and only if
\beq \begin{split} \label{eq:pre_optimized}
	\wt{f}(x) &= c_0 + 2C \left( \frac{2\SNR}{\sqrt{\rat}(w_4-1)} T_1 \left( \frac{x}{2} \right) + \sum_{\ell=2}^{\infty} \frac{1}{\ell}\left(\frac{\SNR}{\sqrt{\rat}}\right)^\ell T_{\ell} \left( \frac{x}{2} \right) \right) \\
	&=  c_0 + 2C \left( \frac{\SNR}{\sqrt{\rat}}\left(\frac{2}{w_4-1}-1\right) T_1 \left( \frac{x}{2} \right) + \sum_{\ell=1}^{\infty} \frac{1}{\ell}\left(\frac{\SNR}{\sqrt{\rat}}\right)^\ell T_{\ell} \left( \frac{x}{2} \right) \right)
\end{split} \eeq
for some constant $c_0$. We notice that the following identity holds for the Chebyshev polynomials:
\beq
\sum_{\ell=1}^{\infty} \frac{t^{\ell}}{\ell} T_{\ell} \left( x \right) = \log \left( \frac{1}{\sqrt{1-2tx+t^2}} \right).
\eeq
(See, e.g., (18.12.9) of \cite{Handbook}.) Since $T_1(x) = x$, we find that \eqref{eq:pre_optimized} is equivalent to
\beq \begin{split}\label{eq:optimized_phi}
	\wt{f}(x) &= c_0 + C \frac{\SNR}{\sqrt{\rat}}\left(\frac{2}{w_4-1}-1\right)x  - C\log\left(\frac{\rat-\SNR\sqrt{\rat}x+w^2}{\rat}\right),
\end{split}\eeq
or
\beq\begin{split}
	f(x) =c_0 + \frac{C \SNR}{\rat}\left(\frac{2}{w_4-1}-1\right)x  -\frac{C \SNR (1+\rat)}{\rat}\left(\frac{2}{w_4-1}-1\right)
	- C\log\left[\frac{\SNR}{\rat}\left(\left(1+\frac{\rat}{\SNR}\right)(1+\SNR)-x\right)\right].
\end{split}\eeq
This concludes the proof of Theorem \ref{thm:CLT} with an optimal function
\beq\label{eq:phi}
\phi_\SNR(x)=\wt \phi_\SNR(\varphi(x))
\eeq 
where
\beq
\wt \phi_w(x)=c_0 +  \frac{\SNR}{\sqrt{\rat}}\left(\frac{2}{w_4-1}-1\right) x  - \log\left(\frac{\rat-\SNR\sqrt{\rat}x+ \SNR^2}{\rat}\right).
\eeq
Choosing
\[
c_0 = \frac{\SNR(1+\rat)}{\rat} \left(\frac{2}{w_4-1}-1\right) + \log (\SNR/\rat),
\]
we get \eqref{eq:optimal_f}.

Next, we prove a lemma for the test statistic in Theorem \ref{thm:test}.
\begin{lem} \label{lem:L_t}
	Let
	\beq
	L_{\SNR} = \sum_{i=1}^M \phi_\SNR(\mu_i) - M \int_{d_-}^{d_+} \phi_{\SNR}(x) \,\dd \mu_{MP}(x)\,,
	\eeq
	where $\phi_w$ is defined as in \eqref{eq:optimal_f}.
	Then
	\beq \begin{split}	
		L_{\SNR} &= -\log \det \left( \left(1+\frac{\rat}{\SNR} \right)(1+\SNR)I - YY^T \right) + \frac{\SNR}{\rat} \left( \frac{2}{w_4-1} - 1 \right) (\Tr YY^T-M) \\
		&\quad + M \left[\frac{\SNR}{\rat} - \log\left(\frac{\SNR}{\rat}\right) -\frac{1-\rat}{\rat}\log(1+\SNR) \right].
	\end{split}	\eeq
\end{lem}

\begin{proof}
	It is straightforward to see that
	\beq
	\sum_{i=1}^M \phi_\SNR(\mu_i) = -\log \det \left( \left(1+\frac{\rat}{\SNR}\right)(1+\SNR)I - YY^T \right) + \frac{\SNR}{\rat} \left( \frac{2}{w_4-1} - 1 \right) \Tr YY^T
	\eeq
	From the well-known formula for the Stieltjes transform of the Marchenko--Pastur law,
	\[
	\int_{d_-}^{d_+} \frac{1}{x-z} \,\dd \mu_{MP}(x) = \frac{1-\rat-z + \sqrt{(z-d_-)(z-d_+)}}{2\rat z}.
	\]
	Integrating it over $z$ and putting $z = (1+\frac{\rat}{\SNR})(1+\SNR)$, we find that
	\[ \begin{split}
	\int_{d_-}^{d_+} \log \left( (1+\frac{\rat}{\SNR})(1+\SNR) - x \right) \dd \mu_{MP}(x) = \frac{\SNR}{\rat} - \log\left(\frac{\SNR}{\rat}\right) -\frac{1-\rat}{\rat}\log(1+\SNR)
	\end{split} \]
	Finally, it is elementary to check that
	\beq
	\int_{d_-}^{d_+} x\;\frac{\sqrt{(x-d_-)(d_+ -x)}}{2\pi \rat x}\dd x= 1.
	\eeq
	This proves the desired lemma.
\end{proof}

Lastly, we prove a lemma for the mean and the variance of the test statistic.
\begin{lem} \label{lem:m_M}
	Let
	\beq
	m_Y(\phi_\SNR)|_{\bsH_0} = \frac{\wt{\phi}_\SNR(2) + \wt{\phi}_\SNR (-2)}{4} -\frac{1}{2} \tau_0(\wt{\phi}_\SNR) - (w_4 -3) \tau_2(\wt{\phi}_\SNR)
	\eeq
	and
	\beq
	m_{Y}(\phi_\SNR)|_{\bsH_1} = \frac{\wt{\phi}_\SNR(2) + \wt{\phi}_\SNR(-2)}{4}  -\frac{1}{2} \tau_0(\wt{\phi}_\SNR) - (w_4 -3) \tau_2(\wt{\phi}_\SNR)  + \sum_{\ell=1}^{\infty} \left(\frac{\SNR}{\sqrt{\rat}}\right)^\ell \tau_{\ell}(\wt{\phi}_\SNR)
	\eeq
	where $\phi_\SNR$ is defined as in \eqref{eq:optimal_f}.
	Then,
	\beq
	m_Y(\phi_\SNR)|_{\bsH_0} = -\frac1{2} \log\left(1-\frac{\SNR^2}{\rat}\right)  -\frac{\SNR^2}{2\rat}(w_4-3)
	\eeq
	and
	\beq
	m_Y(\phi_\SNR)|_{\bsH_1} = m_Y(\phi_\SNR)|_{\bsH_0}  -\log\left(1-\frac{\SNR^2}{\rat}\right) +\frac{\SNR^2}{\rat}\left(\frac2{w_4-1}-1\right).
	\eeq
	In particular, $m_Y(\phi_\SNR)|_{\bsH_0} < m_Y(\phi_\SNR)|_{\bsH_1}$ if $\SNR^2/\rat\in (0, 1)$.
\end{lem}

\begin{proof}
	We first notice that $\wt{\phi}_\SNR$ is the function $\wt{f}$ in \eqref{eq:optimized_phi} with $C=1$ and $c_0 =\tau_0(\wt{\phi}_\SNR) $. Thus, from \eqref{eq:tau_ell},
	\beq \label{eq:tau_ell_phi}
	\tau_1(\wt{\phi}_\SNR) = \frac{2\SNR}{\sqrt{\rat}(w_4-1)}, \quad  \tau_{\ell}(\wt{\phi}_\SNR) = \frac{1}{\ell}\left(\frac{\SNR}{\sqrt{\rat}}\right)^\ell \qquad (\ell = 2, 3, 4, \dots).
	\eeq
	Since
	\beq \begin{split}
		\wt{\phi}_\SNR(2) + \wt{\phi}_\SNR(-2) &= \displaystyle\log\left(\frac{\rat}{\rat-2\SNR\sqrt{\rat}+\SNR^2}\right)+\log\left(\frac{\rat}{\rat+2\SNR\sqrt{\rat}+\SNR^2}\right) + 2c_0 \\
		&= -2 \log\left(1-\frac{\SNR^2}{\rat}\right) + 2c_0,
	\end{split} \eeq
	we find that
	\beq \begin{split}
		m_Y(\phi_w)|_{\bsH_0} =  -\frac1{2} \log\left(1-\frac{\SNR^2}{\rat}\right)  -\frac{\SNR^2}{2\rat}(w_4-3).
	\end{split} \eeq
	Similarly, we also get
	\beq \begin{split}
		m_Y(\phi_\SNR)|_{\bsH_1} &=m_Y(\phi_\SNR)|_{\bsH_0} -\log\left(1-\frac{\SNR^2}{\rat}\right) +\frac{\SNR^2}{\rat}\left(\frac2{w_4-1}-1\right).
	\end{split} \eeq
	Finally, it is obvious $m_Y(\phi_\SNR)|_{\bsH_0} <m_Y(\phi_\SNR)|_{\bsH_1}$ if $\SNR^2/\rat \in (0, 1)$, since $\tau_{\ell}(\wt \phi_\SNR) > 0$ for all $\ell = 1, 2, \dots$.
\end{proof}

\begin{rem} \label{rem:m_M}
	For any $\SNR$, it can be easily checked from \eqref{eq:tau_ell_phi} that
	\beq
	V_{Y}(\phi_\SNR)|_{\bsH_1} = V_Y(\phi_\SNR)|_{\bsH_0} = -2\log\left(1-\frac{\SNR^2}{\rat}\right) +\frac{2\SNR^2}{\rat}\left(\frac2{w_4-1}-1\right),
	\eeq
\end{rem}


\end{document}